\documentclass[11pt,a4paper]{article}
\usepackage{amssymb,amsmath,amsthm,epsfig,verbatim,pstricks,url}
\usepackage{enumitem}
\setlist[itemize]{itemsep=2pt, topsep=0pt}
\usepackage{ifpdf}  
\newtheorem{theorem}{\sc Theorem.}[section]
\newtheorem{lemma}[theorem]{\sc Lemma.}
\newtheorem{remark}[theorem]{\sc Remark.}

\renewcommand{\theequation}{\arabic{section}.\arabic{equation}}
\newenvironment{AMS}%
{{\upshape\bfseries AMS subject classifications. }\ignorespaces}{}
\newenvironment{keywords}{{\upshape\bfseries Key words. }\ignorespaces}{}

\newcommand{\bRplus}{{\mathbb R}_{>0}}
\newcommand{\bRgeq}{{\mathbb R}_{\geq 0}}

\newcommand{\RZ}{{\mathbb R} \slash {\mathbb Z}}
\newcommand{\bC}{{\mathbb C}}
\newcommand{\bR}{{\mathbb R}}
\newcommand{\bS}{{\mathbb S}}
\newcommand{\bN}{{\mathbb N}}
\newcommand{\bZ}{\mathbb{Z}}

\newcommand{\ratio}{{\mathfrak r}}
\newcommand{\dH}[1]{\;{\rm d}{\mathcal{H}}^{#1}} 
\newcommand{\drho}{\;{\rm d}\rho}
\newcommand{\dt}{\;{\rm d}t}
\newcommand{\ds}{\;{\rm d}s}
\newcommand{\Vh}{\underline{V}^h}

\newcommand{\Id}{{\rm Id}}

\newcommand{\dd}[1]{\frac{\rm d}{{\rm d}#1}}
\newcommand{\ddt}{\dd{t}}
\newcommand{\altmu}{{\mathfrak m}}

\def\epsilon{\varepsilon} 

\def\hat{\widehat}

\DeclareMathOperator*{\esssup}{ess\,sup}
\allowdisplaybreaks

\begin{document}
\title{
Discrete anisotropic curve shortening flow in higher codimension}
\author{Klaus Deckelnick\footnotemark[2]\ \and 
        Robert N\"urnberg\footnotemark[3]}

\renewcommand{\thefootnote}{\fnsymbol{footnote}}
\footnotetext[2]{Institut f\"ur Analysis und Numerik,
Otto-von-Guericke-Universit\"at Magdeburg, 39106 Magdeburg, Germany \\
{\tt klaus.deckelnick@ovgu.de}}
\footnotetext[3]{Dipartimento di Mathematica, Universit\`a di Trento,
38123 Trento, Italy \\ {\tt robert.nurnberg@unitn.it}}

\date{}

\maketitle

\begin{abstract}
We introduce a novel formulation for the evolution of parametric curves by
anisotropic curve shortening flow in $\bR^d$, $d\geq2$. The reformulation
hinges on a suitable manipulation of the parameterization's tangential
velocity, leading to a strictly parabolic differential equation.
Moreover, the derived equation is in divergence form, giving rise to a natural
variational numerical method.
For a fully discrete finite element approximation based on piecewise linear
elements we prove optimal error estimates. 
Numerical simulations confirm the theoretical results and demonstrate the
practicality of the method.
\end{abstract} 



\renewcommand{\thefootnote}{\arabic{footnote}}

\setcounter{equation}{0}
\section{Introduction} \label{sec:intro}

The aim of this paper is to define and analyse a finite element approximation 
for the evolution of a curve in $\bR^d$, for an arbitrary $d\geq2$, 
by anisotropic curve shortening flow. The evolution law we consider is a natural  gradient flow for the following anisotropic energy 
\begin{equation} \label{eq:Ephi}
E_\phi(\Gamma)= \int_{\Gamma} \phi(\tau) \dH1, 
\end{equation}
where $\tau$ denotes the unit tangent of $\Gamma$ and $\phi$ is a given, $1$--homogeneous energy density, 
 cf.\ \cite{Pozzi07}. 
For the case $d=2$, i.e.\ curves in the plane, the energy \eqref{eq:Ephi} can
play the role of an interfacial energy, e.g.\ in materials science
\cite{TaylorCH92,Gurtin93}. 
A possible application in differential
geometry is curve shortening flow in a Riemannian manifold.
In this case the energy density in \eqref{eq:Ephi} is required to have 
a spatial dependence, see \cite{BellettiniP96,MaC07}.
For more details on anisotropic surface energies we refer to
\cite{DeckelnickDE05,Giga06} and the references therein.
We shall see in Section~\ref{sec:weak} that a family of curves  $(\Gamma(t))_{t \in (0,T)}$ evolves according to anisotropic curve
 shortening flow provided that
 \begin{equation} \label{eq:anisolaw}
 \mathcal V_{nor} = \varkappa_\phi \quad \mbox{ on } \Gamma(t),
 \end{equation}
where $\mathcal V_{nor}$ denotes the vector of normal velocities, i.e.\
$\mathcal V_{nor} \cdot \tau = 0$, and $\varkappa_\phi$ denotes the 
anisotropic curvature vector. 
In the isotropic case, $\phi(p)=|p|$, we have that 
$\varkappa_\phi=\varkappa=\tau_s$,
with $\cdot_s$ denoting differentiation with respect to arclength, and 
\eqref{eq:anisolaw} is just the usual curve shortening flow. 
If $\Gamma(t)= x(I,t)$ is described by a parameterization
$x:I \times [0,T] \to \bR^d$, where $I=\RZ$ is the periodic interval $[0,1]$,
then \eqref{eq:anisolaw} can be equivalently formulated as
\begin{equation} \label{eq:Pxtflow}
(\Id - \tau \otimes \tau) x_t = \varkappa_\phi,
\end{equation}
where here and throughout we use a slight abuse of notation, in that we
do not distinguish between geometric quantities of the curve being defined
on $\Gamma(t)$ or in $I \times [0,T]$. We observe that \eqref{eq:Pxtflow} 
does not prescribe the tangential velocity. 
\\
{From} a geometric point of view it is natural to consider the system
\begin{equation}\label{eq:anicsf}
x_t = \varkappa_\phi, 
\end{equation}
see e.g.\ \cite{Dziuk94,AmbrosioS96} for the isotropic setting.
Generalizing earlier work \cite{Dziuk94,Dziuk99} by Dziuk, Pozzi \cite{Pozzi07} introduced and analysed a semidiscrete finite element method and obtained
an $\mathcal O(h)$ error bound in $L^\infty(H^1)$ for the position vector. Here, the fact that the anisotropic curvature vector is invariant with respect to reparameterization 
leads  to a degeneracy of the elliptic part on the right hand side of \eqref{eq:anicsf}, thus complicating the error analysis. 
The authors in \cite{curves3d} suggested a discretization of \eqref{eq:Pxtflow} that introduces a tangential velocity at the discrete level which
leads to a nice distribution of vertices in practice. \\
In \cite{finsler}, the present authors considered anisotropic curve shortening flow in the planar case $d=2$. 
The crucial idea from \cite{finsler} was to define positive definite
matrices $H(x_\rho) \in \bR^{2 \times 2}$ such that any solution to the PDE
\begin{equation} \label{eq:finsler}
H(x_\rho) x_t = [\phi(x_\rho) \phi'(x_\rho)]_\rho
\end{equation}
parameterizes anisotropic curve shortening flow. We would like to stress
three appealing aspects of \eqref{eq:finsler}. Firstly, it can be shown that
this PDE is strictly parabolic in the sense of Petrovsky, see
\cite[Lemma~3.3]{finsler}, so that \eqref{eq:finsler} can be interpreted as a kind of DeTurck trick, \cite{ElliottF17}. 
Secondly, the PDE is in divergence form, 
making its variational numerical approximation straightforward, and for a semidiscrete finite element approximation an $\mathcal O(h)$ error bound in $L^\infty(H^1)$ for the
position vector is obtained. 
And thirdly, the tangential velocities induced by \eqref{eq:finsler} lead to
well distributed vertices in practice.
It is the aim of this paper to generalise these results, in that we
\begin{itemize}
\item propose an analogue of \eqref{eq:finsler} for curves in higher codimension;
\item provide an error analysis for a {\em fully} discrete finite element scheme;
\item prove an optimal $L^\infty(L^2)$ error bound.
\end{itemize}
Let us emphasize that all three of the above results are new in the literature.
We note that in \cite{finsler} the anisotropic energy was allowed to depend on
space. We expect that it is possible to treat such a spatial dependence also 
in higher codimensions, allowing to study, for example, curve shortening flow 
in higher dimensional Riemannian manifolds. 
However, in this paper we restrict our attention to 
the simpler energy \eqref{eq:Ephi} in order not to overburden the presentation.

Let us mention
related work on the numerical approximation of (anisotropic) curve shortening
flow in higher codimension. The seminal papers \cite{Dziuk94,DeckelnickD95} 
considered the isotropic case and proved optimal $H^1$-error bounds for
semidiscrete approximations based on piecewise linear elements. 
In the planar case, error estimates for
a fully discrete variant of \cite{Dziuk94}  have recently been obtained in \cite{Li20} and \cite{YeC21}.
The method and numerical analysis in 
\cite{Dziuk94} were generalized to the anisotropic setting in \cite{Pozzi07}.
Both \cite{Dziuk94} and \cite{Pozzi07} may suffer from coalescence of mesh
points in practice, since the discrete curves are only updated in normal
direction.
The so-called BGN schemes, on the other hand, are characterized by an implicit
tangential motion that leads to a nice distribution of vertices, 
\cite{bgnreview}. Their application to (anisotropic) curve shortening flow in 
higher codimensions has been considered in \cite{curves3d}, with an error
analysis for these schemes still lacking.
A semi-Lagrangian scheme in the context of level-set methods was considered in
\cite{CarliniFF07}. A finite volume scheme for possibly interacting
curves driven by curvature forces in $\bR^3$ was introduced in
\cite{BenesKS22}. 
Finally, optimal error estimates for semi- and fully discrete approximations
of a system of PDEs for $\varkappa$ and $\tau$
describing isotropic curve shortening in higher codimension have been obtained
in \cite{BinzK21preprint}.

The remainder of the paper is organised as follows. 
In Section~\ref{sec:weak}, we give a rigorous statement of the 
partial differential equation we wish to study, together with a
derivation of the required matrices $H$.
In Section~\ref{sec:fea} we state a natural weak formulation and
introduce our fully discrete finite element
approximation. We also prove an unconditional stability result for the scheme. 
Section~\ref{sec:proof} is devoted to the proof of our main error 
estimates, which include an $\mathcal O(h+\Delta t)$ bound for a 
discrete $H^1$--norm, and an $\mathcal O(h^2+\Delta t)$ $L^2$--error bound.
Finally, in Section~\ref{sec:nr} we present the results of 
some numerical simulations.

We end this section with a few comments about notation. 
We adopt the standard notation for Sobolev spaces, denoting the norm of
$W^{\ell,p}(I)$ ($\ell \in \bN_0$, $p \in [1, \infty]$)
by $\|\cdot \|_{\ell,p}$ and the 
semi-norm by $|\cdot |_{\ell,p}$. For
$p=2$, $W^{\ell,2}(I)$ will be denoted by
$H^{\ell}(I)$ with the associated norm and semi-norm written as,
respectively, $\|\cdot\|_\ell$ and $|\cdot|_\ell$.
The above are naturally extended to vector functions, and we will write 
$[W^{\ell,p}(I)]^d$ for a vector function with $d$ components.
In addition, we adopt the standard notation $W^{\ell,p}(a,b;B)$
($\ell \in \bN$, $p \in [1, \infty]$, $(a,b)$ an interval in $\bR$, 
$B$ a Banach space) for time dependent spaces
with norm $\|\cdot\|_{W^{\ell,p}(a,b;B)}$.
Once again, we write $H^{\ell}(a,b;B)$ if $p=2$. 
In addition, throughout $C$ denotes a generic positive constant independent of 
the mesh parameter $h$ and the time step size $\Delta t$.
At times $\epsilon$ will play the role of a (small)
positive parameter, with $C_\epsilon>0$ depending on $\epsilon$, but
independent of $h$ and $\Delta t$.
Finally, in this paper we make use of the Einstein summation convention.

\setcounter{equation}{0} 
\section{Mathematical formulation} \label{sec:weak}

\subsection{Anisotropic curve shortening flow} \label{sec:acsf}

Let $\phi \in C^0(\bR^d, \bRgeq) \cap C^4(\bR^d \setminus \{0\}, \bRplus)$,
as well as
\begin{equation} \label{eq:gamma1} 
\phi(\lambda p) =  |\lambda|  \phi(p) \qquad \forall\
p \in \bR^d,\ \lambda \in \bR .
\end{equation}
It is not difficult to verify that \eqref{eq:gamma1} implies that
\begin{equation} \label{eq:phidd}
\phi'(\lambda p) = \frac{\lambda}{|\lambda|} \phi'(p), 
\quad \phi'(p)\cdot p = \phi(p) 
\quad\mbox{and}\quad
\phi''(p) p  = 0 \qquad \forall\ p \in \bR^d\setminus\{0\},\  
\lambda \in \bR \setminus \{0\},
\end{equation}
where $\phi'$ and $\phi''$ denote the gradient and the Hessian 
of $\phi$, respectively.
In addition, we assume that $\phi$ is strictly convex in the sense that
\begin{equation} \label{eq:sconv}
\phi''(p) q \cdot q >0 \qquad \forall\ p,q  \in \bR^d \mbox{ with } |p|=|q| =1, p \cdot q =0.
\end{equation}
Let us consider the anisotropic energy $E_\phi$ defined in terms of $\phi$ via 
\eqref{eq:Ephi} and assume that $\Gamma=\{ x(\rho) : \rho \in I \}$. 
In view of \eqref{eq:gamma1} we have
\begin{displaymath}
E_\phi(\Gamma)= \int_I \phi(\frac{x_\rho}{| x_\rho |}) | x_\rho | \drho = \int_I \phi(x_\rho) \drho.
\end{displaymath}
Using \eqref{eq:phidd} and integration by parts, the first variation of $E_\phi$ in the direction of a vectorfield $V$ defined on $\Gamma$ can be derived as 
\begin{align*}
[\delta E_\phi(\Gamma)] (V) & =
  \dd{\epsilon} \int_I \phi(x_\rho + \epsilon (V\circ x)_\rho) \drho \!\mid_{\epsilon=0} = \int_I \phi'(x_\rho) \cdot (V \circ x)_\rho \drho \\
&= \int_\Gamma \phi'(\tau) \cdot V_s \dH1 = - \int_\Gamma \phi''(\tau) \tau_s \cdot V \dH1
= - \int_{\Gamma} \phi''(\tau) \varkappa \cdot V \dH1,
\end{align*}
where $\varkappa=\tau_s$ denotes the (principal) curvature vector of $\Gamma$. 
The quantity 
\begin{equation} \label{eq:vecvarkappaphi}
{\varkappa}_\phi:
= \phi''(\tau)\varkappa 
\end{equation}
can then be viewed as an anisotropic curvature vector. Note that $\varkappa_\phi \cdot \tau = \varkappa \cdot \phi''(\tau) \tau =0$, so that $\varkappa_\phi \in \{\tau\}^\perp$.  In order to
define a gradient flow for $E_\phi$ we observe that only the normal part of a vectorfield $V$ will contribute to a change in the shape of $\Gamma$. This motivates to consider the gradient
flow with respect to the inner product
\begin{displaymath}
(V,W)_{nor}= \int_{\Gamma} P V \cdot PW \dH1,
\end{displaymath}
see e.g.\ \cite{curves3d}, 
where $P= \Id - \tau \otimes \tau$ denotes the projection onto the normal part of $\Gamma$. A family of curves $(\Gamma(t))_{t \in (0,T)}$ in $\bR^d$ then evolves according to anisotropic
curve shortening flow provided that
\begin{equation} \label{eq:mucsf}
\altmu(\tau) \mathcal{V}_{nor} = \varkappa_\phi,
\end{equation}
where $\altmu : \bS^{d-1} \to \bRplus$ is a given mobility function in order to allow for a more general setting.
For simplicity we extend $\altmu$ 0-homogeneously to $\bR^d \setminus\{0\}$
and require for our analysis that $\altmu\in C^3(\bR^d \setminus\{0\})$.
The case $\altmu(\tau)=1$ is most frequently
treated in the literature, while the choice $\altmu(\tau) = \frac1{\phi(\tau)}$
also has some nice properties and was considered, for example, in
\cite{DeckelnickDE05,Pozzi12,finsler}.  \\
In what follows we assume a parametric description of $\Gamma(t)$, i.e.\ $\Gamma(t)=x(I,t)$, with $I=\bR / \bZ$. Then the unit tangent is
given by  $\tau = \frac{x_\rho}{|x_\rho|}$, while the anisotropic curvature vector is calculated as
\begin{displaymath}
\varkappa_\phi = \phi''(\tau) \varkappa =  \phi''(\tau) \tau_s = [\phi'(\tau)]_s = \frac{1}{| x_\rho|} [\phi'(\tau)]_\rho.
\end{displaymath}
Thus, the family $(\Gamma(t))_{t \in (0,T)}$ evolves according to \eqref{eq:mucsf} provided that 
\begin{equation} \label{eq:iow}
\altmu(\tau) P x_t = \frac{[\phi'(\tau)]_\rho}{|x_\rho|} .
\end{equation}

\subsection{DeTurck's trick for anisotropic curve shortening flow} 
\label{sec:DTT}

In what follows we shall make frequent use of the function
\begin{equation*} 
\Phi(p)=\tfrac12 \phi^2(p).
\end{equation*}
It is not difficult to verify that 
$\Phi \in C^1(\bR^d) \cap C^4(\bR^d \setminus \{0\})$ is convex and that
\begin{equation} \label{eq:Phidd}
\Phi''(\lambda p) = \Phi''(p), \quad
\Phi''(p) p = \Phi'(p) \quad \text{and} \quad \Phi'''(p) 
(p, \cdot,\cdot) = 0
\qquad \forall\ p \in \bR^d\setminus\{0\},\  
\lambda \in \bR \setminus\{0\} ,
\end{equation}
where we think of $\Phi'''(p)$ as a symmetric trilinear form on $\bR^d \times \bR^d \times \bR^d$.
Moreover, we have the following result.
\begin{lemma}
For every compact set $K \subset \bR^d \setminus \{0\}$ there exists $\sigma=\sigma(K)>0$ such that 
\begin{subequations} \label{eq:convex}
\begin{alignat}{2}
\Phi''(p) \xi \cdot \xi & \geq \sigma | \xi |^2 \quad \quad 
&& \forall\ p \in K,\ \xi \in \bR^d, \label{eq:convex3} \\
\Phi(q) - \Phi(p) - \Phi'(p) \cdot (q-p) & \geq \sigma | q-p |^2 \quad 
&& \forall\ p,q \in K \mbox{ with } [p,q] \subset K, \label{eq:convex2} 
\end{alignat}
\end{subequations}
where $[p,q]$ denotes the segment between $p$ and $q$.
\end{lemma}
\begin{proof}  
It is shown in \cite[Remark~1.7.5]{Giga06} that \eqref{eq:sconv} implies that
$\Phi''(p)$ is positive definite for all $p \in \bR^d \setminus \{0\}$.
Hence the estimates \eqref{eq:convex3} and \eqref{eq:convex2}
follow.
\end{proof}

As described in the introduction, our aim is  to construct positive definite matrices 
$H(p) \in \bR^{d \times d}$ such that solutions of the system
\begin{equation} \label{eq:Hxt}
H(x_\rho) x_t = [\Phi'(x_\rho)]_\rho 
\end{equation}
solve \eqref{eq:iow}. 
To this end, for any $p \in \bR^{d} \setminus \{0\}$ we make the ansatz
\begin{equation} \label{eq:newHminusT}
 H^{-1}(p) = \frac1{\alpha(p)} \left( \Id 
+ \mathfrak t \otimes w(p) - w(p) \otimes\mathfrak t\right) ,
\quad \text{where} \quad \mathfrak t = \frac{p}{|p|},
\end{equation}
and where $\alpha(p) \in \bRplus$ and $w(p) \in \bR^d$ have to be chosen
appropriately.
Here the two terms $\Id + \mathfrak t \otimes w$ in the definition of 
$H^{-1}$ are guided by the insight that tangential changes to $x_t$ do
not change the parameterized flow.
Adding the third term $- w\otimes\mathfrak t$ in 
\eqref{eq:newHminusT} then ensures positive definiteness. 
Remarkably, the simple scaling
factor $1/\alpha$ allows enough freedom to guarantee that in the normal
directions the correct flow is obtained, i.e.\ that \eqref{eq:iow} holds.

We will now derive values for $\alpha(p)$ and $w(p)$ so that \eqref{eq:iow} 
holds. If we assume that \eqref{eq:Hxt} holds, 
then the ansatz \eqref{eq:newHminusT} yields
\begin{align} \label{eq:ansatz1}
P x_t & = P H^{-1}(x_\rho) [\Phi'(x_\rho)]_\rho 
= \frac1{\alpha} (\Id - \tau \otimes \tau) \left( 
\Id + \tau \otimes w - w \otimes \tau \right)
[\Phi'(x_\rho)]_\rho \nonumber \\ &
= \frac1{\alpha} \left(
\Id + \tau \otimes w - w \otimes \tau
- \tau \otimes \tau - \tau \otimes w + (\tau \cdot w) \tau \otimes \tau
\right) [\Phi'(x_\rho)]_\rho \nonumber \\ &
= \frac1{\alpha} \left(
P - w \otimes \tau + (\tau \cdot w) \tau \otimes \tau \right) 
[\Phi'(x_\rho)]_\rho.
\end{align}
Direct calculation shows that
\begin{equation} \label{eq:newHxt1}
[\Phi'(x_\rho)]_\rho 
= [\phi(x_\rho) \phi'(x_\rho)]_\rho 
= [\phi(x_\rho)]_\rho \phi'(x_\rho) + \phi(x_\rho) [ \phi'(x_\rho)]_\rho.
\end{equation}
Combining \eqref{eq:ansatz1} and \eqref{eq:newHxt1}, on noting 
$[\phi'(x_\rho)]_\rho \cdot \tau = 0$, recall \eqref{eq:vecvarkappaphi}, 
and \eqref{eq:phidd}, yields that
\begin{align*}
P x_t & = \frac1{\alpha} \left(
[\phi(x_\rho)]_\rho P \phi'(x_\rho) + \phi(x_\rho) P [ \phi'(x_\rho)]_\rho
\right. \nonumber \\ & \left. \qquad
 - [\phi(x_\rho)]_\rho (\phi'(x_\rho)\cdot \tau) w 
+ [\phi(x_\rho)]_\rho (\tau \cdot w) (\phi'(x_\rho) \cdot \tau) \tau \right) 
\nonumber \\ &
= \frac1{\alpha} \phi(x_\rho) [ \phi'(x_\rho)]_\rho
+ \frac{[\phi(x_\rho)]_\rho}{\alpha}
\left(
 P \phi'(x_\rho) - \phi(\tau) w + (\tau \cdot w) \phi(\tau) \tau
\right) \nonumber \\ &
= \frac1{\alpha} \phi(x_\rho) [ \phi'(x_\rho)]_\rho, 
\end{align*}
provided that
\begin{equation*} 
P \phi'(x_\rho) - \phi(\tau) w + (\tau \cdot w) \phi(\tau) \tau = 0,
\end{equation*}
which can be achieved by setting $w = \frac1{\phi(\tau)} P \phi'(x_\rho)
= \frac{|x_\rho|}{\phi(x_\rho)} (\Id - \frac{x_\rho\otimes x_\rho}{|x_\rho|^2})
\phi'(x_\rho)$.
Hence with this choice \eqref{eq:iow} will be satisfied if we let
\begin{equation*} 
\alpha = \altmu(\tau)|x_\rho|\phi(x_\rho).
\end{equation*}

Before we summarize our results, we state some properties of $H$
that immediately follow from the ansatz \eqref{eq:newHminusT}.

\begin{lemma} \label{lem:posdef}
For $\mathfrak t \in \bS^{d-1}$, $w \in \{\mathfrak t\}^\perp$ 
and $\alpha \in \bRplus$ let
\begin{equation} \label{eq:lemHinv}
H^{-1} = \frac1{\alpha} \left( \Id + \mathfrak t \otimes w 
- w \otimes \mathfrak t \right).
\end{equation}
Then it holds that
\begin{equation} \label{eq:lemH}
H 
= \alpha \left[\Id + \frac{w\otimes (\mathfrak t - w ) 
- \mathfrak t \otimes (w + |w|^2 \mathfrak t)}{1 +  |w|^2} \right] 
.
\end{equation}
Moreover, the matrix $H$ is positive definite and satisfies
\begin{equation} \label{eq:newHposdef}
H z \cdot z \geq \frac{\alpha}{1+|w|^2}|z|^2 
\qquad \forall\ z \in \bR^d.
\end{equation}
\end{lemma}
\begin{proof}
Direct calculation shows that
\begin{align*}
H^{-1} H & = 
\frac1{\alpha} 
\left( H + \mathfrak t \otimes w H - w\otimes\mathfrak t H 
\right) \nonumber \\ & 
= \frac1\alpha H + \frac{\mathfrak t \otimes [ (1 +  |w|^2) w
+ |w|^2 (\mathfrak t - w) ] - w \otimes [ (1 +  |w|^2) \mathfrak t
- w - |w|^2 \mathfrak t]}{1 +  |w|^2} \nonumber \\ & 
= \frac1\alpha H + \frac{\mathfrak t \otimes [ w
+ |w|^2 \mathfrak t ] - w \otimes [ \mathfrak t - w ]}
{1 +  |w|^2}
 = \Id, 
\end{align*}
which proves that \eqref{eq:lemH} is indeed the inverse of \eqref{eq:lemHinv}. 
Moreover, it holds for $w \neq 0$ that 
\begin{align*}
H z \cdot z & = \alpha |z|^2
- \frac{\alpha}{1+|w|^2} \left( 
(z \cdot w)^2 + |w|^2 (z \cdot \mathfrak t)^2 \right) = \alpha |z|^2
- \frac{\alpha |w|^2}{1+|w|^2} \left( 
(z \cdot \frac{w}{|w|})^2 + (z \cdot \mathfrak t)^2 \right)   \nonumber \\ &
\geq 
\alpha |z|^2
- \frac{\alpha |w|^2}{1+|w|^2} |z|^2
= \frac{\alpha}{1+|w|^2} |z|^2 , 
\end{align*}
where we also used that $w \cdot \mathfrak t=0$.
This proves the desired result \eqref{eq:newHposdef}, since it holds trivially
for $w=0$. 
\end{proof}

In summary, we have shown the following result.

\begin{theorem} 
Let $\Phi(p) = \tfrac12 \phi^2(p)$ and define for
$p \in\bR^d\setminus \{0\}$
\begin{equation} \label{eq:Hfinal}
H(p)
= \alpha(p) \left[\Id + \frac{w(p)\otimes ( \frac{p}{|p|} - w(p) ) 
- \frac{p}{|p|} \otimes (w(p) + |w(p)|^2 \frac{p}{|p|})}{1 +  |w(p)|^2} 
\right] ,
\end{equation}
with
\begin{equation} \label{eq:alphawfinal}
\alpha(p) = \altmu(p)|p|\phi(p) > 0 \quad\text{and}\quad
w(p) = \frac{|p|}{\phi(p)} \left(\Id - \frac{p\otimes p}{|p|^2}\right) \phi'(p)
= \frac{|p|}{\phi(p)}\phi'(p) - \frac{p}{|p|} .
\end{equation}
If $x : I \times [0,T] \to \bR^d$ satisfies \eqref{eq:Hxt},
then $x$ is a solution to anisotropic curve shortening flow, \eqref{eq:iow}. 
Furthermore, $p \mapsto H(p)$ belongs to $C^3(\bR^d \setminus \{0\},
\bR^{d\times d})$ and $H(p)$ is positive definite with
\begin{equation} \label{eq:Hposdef}
H(p) z \cdot z \geq \frac{\alpha(p)}{1 + |w(p)|^2} |z|^2 \qquad \forall\
p \in \bR^d \setminus \{0\},\ z \in \bR^d.
\end{equation}
\end{theorem}

In addition, we can show that \eqref{eq:Hxt} is parabolic.
\begin{lemma} \label{lem:petrovsky}
The system \eqref{eq:Hxt} is parabolic in the sense of Petrovsky.
\end{lemma}
\vspace{-4mm}
\begin{proof}
On inverting the matrix $H(x_\rho)$ we may write \eqref{eq:Hxt} in the form
\begin{equation*} 
 x_t = H^{-1}(x_\rho) \Phi''(x_\rho) x_{\rho \rho} 
\quad \text{in } I \times (0,T].
\end{equation*}
Hence, by definition we need to show that the eigenvalues
of $H^{-1}(p) \Phi''(p)$ have positive real parts for every 
$p \in \bR^d \setminus \{0\}$, 
see e.g.\ \cite[Definition~1.2]{EidelmanIK04}. Let us fix $p \in \bR^d \setminus \{0\}$.
A straightforward extension of the proof of \eqref{eq:newHposdef} shows that 
\begin{displaymath}
\mbox{Re} (H(p) z \cdot \bar z) \geq \frac{\alpha(p)}{1+|w(p)|^2}|z|^2 
\qquad \forall\ z \in \bC^d.
\end{displaymath}
Let $(\lambda,z) \in \bC \times \bC^d \setminus \{0\}$ be an eigenpair of 
$H^{-1}(p) \Phi''(p)$, i.e.\ $\Phi''(p)z=\lambda H(p)z$. Then we have
\begin{displaymath}
\mbox{Re} \lambda= \mbox{Re} \left( \frac{\Phi''(p) z \cdot \bar z}{H(p) z \cdot \bar z} \right) = \frac{\Phi''(p) z \cdot \bar z}{| H(p) z \cdot \bar z |^2} 
\mbox{Re} (H(p) z \cdot \bar z) >0,
\end{displaymath}
where we have used that $\Phi''(p)$ is symmetric and positive definite.
\end{proof}

\begin{remark} \label{rem:d2}
a) In the isotropic case $\phi(p) = |p|$ with $\altmu = 1$, we have 
$w(p) = 0$ and $\alpha(p) = |p|^2$, so that $H(p)=|p|^2 \Id$. Therefore
\eqref{eq:Hxt} becomes
$| x_\rho |^2 x_t = x_{\rho \rho}$, which is precisely the equation considered in \cite{DeckelnickD95}. \\[2mm]
b) We would like to compare \eqref{eq:Hfinal} with the 2d analogue from
\cite{finsler}. In fact, for the setting in \cite{finsler} we have
$\altmu(p) = \frac{|p|}{\phi(p)}$ and
$\phi(p)=\gamma(p^\perp)$, where 
$p^\perp= \binom{p_1}{p_2}^\perp = \binom{-p_2}{p_1}$
and $\gamma:\bR^2 \to \bRgeq$ is a normal-dependent
anisotropic density function.
Hence we have from \eqref{eq:alphawfinal} that $\alpha(p) = |p|^2$ and
\[
w= \frac{| p |}{\phi(p)} (\Id - \mathfrak t \otimes \mathfrak t) \phi'(p) 
= \frac{| p |}{\phi(p)} \bigl( \phi'(p) \cdot \mathfrak t^\perp \bigr) \mathfrak t^\perp = - \frac{1}{\gamma(p^\perp)} \bigl( \gamma'(p^\perp)^\perp \cdot p^\perp \bigr) 
\mathfrak t^\perp = - \frac{1}{\gamma(p^\perp)} \bigl( \gamma'(p^\perp) \cdot p \bigr) \mathfrak t^\perp.
\]
On the other hand, \eqref{eq:lemH}, with $\alpha=|p|^2$, can be re-written as
\begin{align} \label{eq:remH2}
H & =  
 |p|^2 \left[ \Id + \frac{w\otimes \mathfrak t - \mathfrak t \otimes w 
- (w \otimes w + |w|^2 \mathfrak t \otimes \mathfrak t)}{1 +  |w|^2} \right] 
\nonumber \\
& = \frac{|p|^2}{1 +  |w|^2} \left [\Id + w\otimes \mathfrak t - \mathfrak t \otimes w
\right] = \frac{|p|^2}{1+|w|^2} \begin{pmatrix} 1 & - w \cdot \mathfrak t^\perp \\ w \cdot \mathfrak t^\perp & 1 \end{pmatrix} , 
\end{align}
where we have used that for $d=2$ it holds that $w \otimes w + |w|^2 \mathfrak t \otimes \mathfrak t = |w|^2 \Id$. Observing that $w \cdot \mathfrak t^\perp = - \frac{1}{\gamma(p^\perp)} \bigl( \gamma'(p^\perp) \cdot p \bigr)$
and
\begin{equation*} 
1 + |w|^2 = 1 + \frac{(\gamma'(p^\perp)  \cdot p)^2}{\gamma^2(p^\perp)} 
= \frac{1}{\gamma^2(p^\perp)} \left[ (\gamma'(p^\perp) \cdot p^\perp)^2 +
(\gamma'(p^\perp) \cdot p)^2 \right]
= \frac{| p |^2 |\gamma'(p^\perp)|^2}{\gamma^2(p^\perp)}
\end{equation*}
we finally obtain from \eqref{eq:remH2} that
\begin{equation*} 
 H = \frac{\gamma(p^\perp)}{|\gamma'(p^\perp)|^2}
\begin{pmatrix}
\gamma(p^\perp) &  \gamma'(p^\perp) \cdot p \\
-\gamma'(p^\perp) \cdot p & \gamma(p^\perp)
\end{pmatrix} ,
\end{equation*}
which is precisely (3.21b) from \cite{finsler}. 
\end{remark}

\setcounter{equation}{0}
\section{Fully discrete finite element approximation} \label{sec:fea}

The weak formulation corresponding to \eqref{eq:Hxt} reads as follows.
Given $x_0 : I \to \bR^d$, find $x:I \times [0,T] \rightarrow \bR^d$ such 
that $x(\cdot,0)= x_0$ and, for $t \in (0,T]$,
\begin{equation} \label{eq:Hxtweak}
\int_I H(x_\rho) x_t \cdot \eta \drho 
+ \int_I \Phi'(x_\rho) \cdot \eta_\rho \drho =0 
\qquad \forall\ \eta \in [H^1(I)]^d.
\end{equation}
On choosing $\eta = x_t$ in \eqref{eq:Hxtweak} we obtain the natural
energy estimate
\begin{equation} \label{eq:stab}
\ddt \int_I \Phi( x_\rho) \drho 
= - \int_I  H(x_\rho) x_t \cdot x_t \drho \leq 0,
\end{equation}
in view of \eqref{eq:Hposdef}.

In order to define our finite element approximation, 
let $0=q_0 < q_1 < \ldots < q_{J-1}< q_J=1$  be a
decomposition of $[0,1]$ into intervals $I_j=[q_{j-1},q_j]$. Let $h_j=q_j - q_{j-1}$ as well as $h=\max_{1 \leq j \leq J} h_j$. We
assume that there exists a positive constant $c$ such that
\begin{equation*} 
h \leq c h_j, \quad 1 \leq j \leq J,
\end{equation*}
so that the resulting family of partitions of $[0,1]$ is quasi-uniform. 
Within $I$ we identify $q_J=1$ with $q_0=0$ and define the finite element 
spaces
\begin{displaymath}
V^h = \{\chi \in C^0(I) : \chi\!\mid_{I_j} \mbox{ is affine},\ j=1,\ldots, J\}
\quad\text{and}\quad \Vh = [V^h]^d.
\end{displaymath}
Let $\{\chi_j\}_{j=1}^J$ denote the standard basis of $V^h$.
For later use, we let $\pi^h:C^0(I)\to V^h$ 
be the standard interpolation operator at the nodes $\{q_j\}_{j=1}^J$,
and we use the same notation for the interpolation of vector-valued
functions. 
It is well-known that for 
$k \in \{ 0,1 \}$, $\ell \in \{ 1,2 \}$ and $p \in [2,\infty]$ it holds that
\begin{subequations}
\begin{alignat}{2}
h^{\frac 1p - \frac 1r} \| \eta_h \|_{0,r} 
+ h | \eta_h |_{1,p} & \leq C \| \eta_h \|_{0,p} 
\qquad && \forall\ \eta_h \in V^h, \qquad r \in [p,\infty], 
\label{eq:inverse} \\
| \eta - \pi^h \eta |_{k,p} & \leq Ch^{\ell-k} | \eta |_{\ell,p} 
\qquad && \forall\ \eta \in W^{\ell,p}(I). \label{eq:estpih} 
\end{alignat}
\end{subequations}

In order to discretize in time, let $t_m=m \Delta t$, $m=0,\ldots,M$, 
with the uniform time step $\Delta t = \frac TM >0$. 
Then our finite element scheme is defined as follows.
Let $x^0_h = Q_h x_0$, where the nonlinear projection $Q_h$
is defined in Lemma~\ref{lem:nonlinearinterpol} below. 
Then, for $m=0,\ldots,M-1$, find $x^{m+1}_h \in \Vh$ such that
\begin{equation}\label{eq:fea}
\frac1{\Delta t} \int_I H(x^m_{h,\rho}) (x^{m+1}_{h} - x^m_h) \cdot \eta_h 
\drho
+ \int_I \Phi'(x^{m+1}_{h,\rho})\cdot \eta_{h,\rho} \drho =0
\qquad \forall\ \eta_h \in \Vh.
\end{equation}
We note that in the case $d=2$, and for the mobility
$\altmu(\tau)=\frac1{\phi(\tau)}$, the scheme \eqref{eq:fea},
with mass lumping used in the first integral, is identical to
the fully discrete approximation (5.4) from \cite{finsler}, recall also
Remark~\ref{rem:d2}.

We begin by stating an unconditional stability result for the scheme
\eqref{eq:fea}. 

\begin{theorem} 
Any solution of \eqref{eq:fea} satisfies the energy estimate 
\begin{equation} \label{eq:fdstab}
\int_I \Phi( x^{m+1}_{h,\rho}) \drho
+ \frac1{\Delta t} \int_I
H(x^m_{h,\rho}) (x^{m+1}_h - x^m_h) \cdot (x^{m+1}_h - x^m_h) \drho
\leq \int_I \Phi( x^m_{h,\rho}) \drho ,
\end{equation}
for $m=0,\ldots,M-1$.
\end{theorem}
\begin{proof}
The convexity of $\Phi$ implies that $\Phi'(p) \cdot (p-q)  \geq \Phi(p) - \Phi(q)$ for all
$ p,q \in \bR^d$ so that
\begin{equation} \label{eq:Phiconvex}
\int_I \Phi'(x^m_{h,\rho}+\eta_{h,\rho}) \cdot \eta_{h,\rho} \drho
\geq \int_I \Phi(x^m_{h,\rho}+\eta_{h,\rho}) - \Phi(x^m_{h,\rho}) \drho,
\end{equation}
for all $\eta_h \in \Vh$. Choosing $\eta_h = x^{m+1}_h - x^m_h$ in
\eqref{eq:fea} and applying \eqref{eq:Phiconvex} yields the bound
\eqref{eq:fdstab}. 
\end{proof}

Observe that \eqref{eq:fdstab} is a fully discrete analogue of \eqref{eq:stab}.
We note that a discrete analogue of $\ddt \int_I \phi(x_\rho) \drho \leq 0$,
in analogy to this property holding for solutions of the continuous problem
\eqref{eq:Hxtweak}, is much harder to prove. 
For the isotropic case such a discrete analogue can be found in 
\cite[Lemma~4.1.3]{BanschDGP23} for the scheme proposed
in \cite{DeckelnickD95} with mass lumping.
However, extending these techniques to the anisotropic
problem studied here appears to be highly nontrivial.
Nevertheless, we remark that in 
all our numerical experiments, both $\tfrac12 \int_I \Phi(x^m_{h,\rho}) \drho$
and $\int_I \phi(x^m_{h,\rho}) \drho$ are monotonically decreasing. 

Our main result is stated in the following theorem. Here, and from now on,
for a function $f \in C([0,T];B)$, with some Banach space $B$, 
we let $f^m= f(t_m)$.

\begin{theorem} \label{thm:main}
Suppose that \eqref{eq:Hxtweak} has a smooth solution 
$x:I \times [0,T] \rightarrow \bR^d$ satisfying 
\begin{equation} \label{eq:xregul}
x \in  W^{1,\infty}(0,T; [W^{2,\infty}(I)]^d), \, x_{tt} \in L^\infty(0,T; [L^\infty(I)]^d) 
\end{equation}
and
\begin{equation} \label{eq:regul}
c_0 \leq | x_\rho | \leq C_0 \quad \mbox{ in } I \times [0,T]
\end{equation}
for some constants $c_0, C_0 \in \bRplus$.
Then there exist $\delta>0$ and $h_0 > 0$ such that
if $0<h \leq h_0$ and $\Delta t \leq \delta h$,
then \eqref{eq:fea} has a unique 
solution $(x^{m}_h)_{m=0,\ldots,M}$, and the following error bounds hold:
\begin{equation} \label{eq:eb}
\max_{m=0,\ldots,M} \| x^m - x^m_h \|_0^2 
\leq C \left( h^4 + (\Delta t)^2 \right), \quad 
\max_{m=0,\ldots,M} | x^m - x^m_h |_1^2 
\leq C \left( h^2 + (\Delta t)^2 \right).
\end{equation}
\end{theorem}

\setcounter{equation}{0}
\section{Proof of Theorem \ref{thm:main}} \label{sec:proof}

As a crucial ingredient of our error analysis we introduce the following 
nonlinear Ritz--type projection.

\begin{lemma} \label{lem:nonlinearinterpol}
Let $y \in [W^{2,\infty}(I)]^d$ with $\Vert y \Vert_{2,\infty} \leq C_1$ and $c_1 \leq | y_\rho | \leq C_1$ in $I$ for 
some $C_1,c_1>0$. Then there exists a unique function $Q_h y \in \Vh$ such that
\begin{equation} \label{eq:defhatx}
\int_I \Phi''(y_\rho) (Q_h y - y)_\rho \cdot \eta_{h,\rho} \drho
+ \int_I (Q_h y - y) \cdot \eta_h \drho = 0
\quad \forall\ \eta_h \in \Vh.
\end{equation}
Furthermore, there exist $h_1>0$ and $C > 0$ depending on $C_1, c_1$ and $\Phi$ such that
\begin{subequations} \label{eq:approxboth}
\begin{align} 
& \| y - Q_h y \|_0 + h | y - Q_h y |_1 \leq C h^2, \label{eq:approx} \\
& | y - Q_h y |_{1,\infty} \leq C h, \quad 
\tfrac12 c_1 \leq | (Q_h y)_\rho | \leq 2 C_1 
\text{ a.e.\ in } I \label{eq:approxinf}
\end{align}
\end{subequations}
for $0<h \leq h_1$. 
Similarly, if $y \in W^{1,\infty}(0,T; [W^{2,\infty}(I)]^d)$ with 
$\esssup_{0 < t < T} \Vert (y,y_t)(\cdot,t) \Vert_{2,\infty} \leq C_1$ and
$| y_\rho | \geq c_1$ a.e.\ in $I \times [0,T]$, 
then there exist $h_1>0$ and  $C$ depending on $C_1, c_1$ and $\Phi$ such that
for all $0<h \leq h_1$
\begin{equation} \label{eq:approxtime}
\| y_t - (Q_h y)_t \|_0 + h | y_t -  (Q_h y)_t |_1 \leq C h^2 \quad\text{and}\quad | (Q_h y)_t |_{1,\infty} \leq C
\quad \text{a.e.\ in}\ [0,T].
\end{equation}
\end{lemma}
\begin{proof} See Appendix~\ref{sec:appA}.
\end{proof}

Let us define
\begin{equation*} 
E^m = \int_I \Phi(x^m_{h,\rho}) - \Phi((Q_h x^m)_\rho) 
- \Phi'( (Q_h x^m)_\rho) \cdot (x^m_h - Q_h x^m)_\rho \drho,
\end{equation*}
as well as
\begin{equation*} 
F^m=  \int_I \left[ ( (Q_h x^m)_i - x^m_i) \xi^m_i - \zeta^m \right] 
\cdot (Q_h x^m - x^m_h)_\rho \drho,
\end{equation*}
where we have used Einstein summation convention and where
\begin{align} \label{eq:defxizeta}
\xi^m_i & = H_{p_{i}}(x^m_\rho) x^{m+1}_t,\ i = 1,\ldots,d, \quad \nonumber \\
\text{and }\ 
\zeta^m & = \zeta(\cdot,t_m) \ \text{ for }\ 
\zeta= \Phi'((Q_h x)_\rho) - \Phi'(x_\rho) - \Phi''(x_\rho)(Q_h x - x)_\rho.
\end{align}
We note in view of \eqref{eq:convex2} that $E^m$ behaves like $|Q_h x^m - x^m_h |_1^2$, while $F^m$ does not have a sign but
will be controlled by $E^m$ and with the help of Lemma~\ref{lem:nonlinearinterpol}.
Our aim is to obtain
the superconvergence bound $\mathcal O(h^4+(\Delta t)^2)$ for $E^m$, which in turn will lead to the optimal $L^2$--error estimate in \eqref{eq:eb}. 
This procedure can be seen as a nonlinear variant of a technique introduced by Wheeler, \cite{Wheeler73}, for the heat equation.

We shall prove Theorem~\ref{thm:main} with the help of an induction argument. 
In particular, we will prove that there exist $h_0 >0$, $0< \delta \leq 1$ and 
$\mu>0$ such that if $0< h \leq h_0$ and $\Delta t \leq \delta h$, 
then for $m \in \{ 0,\ldots,M \}$ the discrete solution $x^m_h$ exists 
and satisfies
\begin{equation}  \label{eq:ind}
| x^m_{h,\rho} - x^m_\rho | \leq \tfrac12 c_0 \quad \text{in}\ I,  
\quad E^m+ F^m \leq \bigl( h^4 + (\Delta t)^2  \bigr) e^{\mu t_m}.
\end{equation}
Since $x^0_h = Q_h x^0$, the assertion 
\eqref{eq:ind} clearly holds for $m=0$ in view of \eqref{eq:regul} and 
\eqref{eq:approxboth}, 
for $h_0$ chosen sufficiently small and for
arbitrary $0 < \delta \leq 1$ and $\mu>0$.
On assuming that \eqref{eq:ind} holds for a fixed 
$m \in \{ 0,\ldots,M-1 \}$, we will now show that it also holds for $m+1$. 
Let us define 
\[
K= \{ z \in \bR^d : \tfrac14 c_0 \leq | z | \leq 4 C_0 \}.
\]
We infer from \eqref{eq:convex}, \eqref{eq:Hposdef} and 
\eqref{eq:alphawfinal} that there exists $\sigma>0$ such that
\eqref{eq:convex} holds for this $K$, as well as
\begin{equation} \label{eq:Hposdef1}
H(p) \xi \cdot \xi  \geq \sigma | \xi |^2 \quad \quad 
 \forall\ p \in K,\ \xi \in \bR^d. 
\end{equation}
Let us abbreviate
\begin{equation} \label{eq:defeh}
e^m_h=Q_h x^m - x^m_h.
\end{equation}
We have for any $z \in [x^m_{h,\rho}(\rho),(Q_h x^m)_\rho(\rho)]$, 
say $z=\lambda x^m_{h,\rho}(\rho)+(1-\lambda) (Q_h x^m)_\rho(\rho)$, that
\begin{equation} \label{eq:z}
| z | \geq | x^m_\rho | - (1-\lambda) | (Q_h x^m)_\rho - x^m_\rho| - \lambda | x^m_{h,\rho}-x^m_\rho | \geq c_0 -Ch -\tfrac12 c_0 \geq \tfrac14 c_0,
\end{equation}
provided that $h_0$ is small enough,
where we have used \eqref{eq:regul}, \eqref{eq:ind} and \eqref{eq:approxinf}.
Deriving an upper bound for $z$ in a similar way one obtains that $[x^m_{h,\rho}(\rho),(Q_h x^m)_\rho(\rho)] \subset K$ for
all $\rho \in I$. Thus we deduce with the help
of \eqref{eq:convex2}, \eqref{eq:approxboth} and \eqref{eq:ind} that
\begin{align*}
\sigma | e^m_h |_1^2 \leq E^m &
\leq E^m + F^m + C \left( \| x^m - Q_h x^m \|_0 
+ | x^m - Q_h x^m |_1 | x^m - Q_h x^m |_{1,\infty} \right) |e^m_h|_1 \\ &
\leq E^m+ F^m + C h^2 | e^m_h |_1 
\leq E^m + F^m + \tfrac12 \sigma | e^m_h |_1^2 + C h^4, 
\end{align*}
and therefore
\begin{equation} \label{eq:h1err}
| e^m_h |_1^2 \leq C\bigl( E^m+ F^m \bigr) + C h^4 \leq 
C \bigl( h^4 + (\Delta t)^2  \bigr) e^{\mu T} 
 \leq C \bigl( h^4  +  \delta^2 h^2  \bigr) e^{\mu T}  \leq C h^2
\end{equation} 
provided that $ h_0^2 e^{\mu T} \leq 1$ and $\delta^2 e^{\mu T} \leq 1$. 
Let us now begin with the induction step. 

\begin{lemma} \label{lem:ex}
There exists $\delta>0$ such that for $\Delta t \leq \delta h$ there exists a unique element $x^{m+1}_h \in \Vh$ satisfying \eqref{eq:fea} as well as
\begin{equation} \label{eq:timedif1}
| x^{m+1}_h - x^m_h |_{1,\infty} \leq C(h^{\frac12}+ \delta^{\frac12}).
\end{equation}
\end{lemma}
\begin{proof} 
First of all, the existence and uniqueness 
of $x^{m+1}_h$ can be obtained as in \cite[Theorem~2.3]{eqdproc},
taking advantage of the inequality
\begin{equation} \label{eq:Phiconvex3}
(\Phi'(q) - \Phi'(p)) \cdot (q-p) \geq \hat{c} | q-p|^2 
\qquad \forall\ p,q \in \bR^d, 
\end{equation}
where $\hat c > 0$ a positive constant. We remark that
the proof of \eqref{eq:Phiconvex3} from
\cite[Lemma~2.2]{eqdproc} can be easily generalized to $d\geq2$.

Next, on choosing $\eta_h=x^{m+1}_h - x^m_h$ in \eqref{eq:fea}, and
using \eqref{eq:Hposdef1}, \eqref{eq:Phiconvex3} together with the fact that
$x^m_{h,\rho} \in K$ we deduce that
\begin{align}\label{eq:timedif}
& 
\frac{\sigma}{\Delta t} \| x^{m+1}_h - x^m_h \|_0^2 + \hat{c} | x^{m+1}_h- x^m_h |_1^2 \nonumber \\ & \
 \leq \frac1{\Delta t}\int_I H(x^m_{h,\rho}) (x^{m+1}_h - x^m_h) \cdot
(x^{m+1}_h - x^m_h) \drho + \int_I \bigl( \Phi'(x^{m+1}_{h,\rho}) - \Phi'(x^m_{h,\rho}) \bigr) \cdot (x^{m+1}_{h,\rho} - x^m_{h,\rho}) \drho \nonumber \\ & \
 = - \int_I \Phi'(x^m_{h,\rho}) 
\cdot (x^{m+1}_{h,\rho}- x^m_{h,\rho}) \drho \nonumber \\ & \
=  \int_I \bigl( \Phi'(x^m_{\rho}) - \Phi'(x^m_{h,\rho}) \bigr) 
\cdot (x^{m+1}_{h,\rho}- x^m_{h,\rho}) \drho
+  \int_I (\Phi'(x^m_\rho))_\rho  \cdot (x^{m+1}_h- x^m_h) \drho =: S_1 + S_2
,
\end{align}
where we have used integration by parts in the last step.
Since $0 \notin [x^m_\rho(\rho),x^m_{h,\rho}(\rho)] \subset K$ for all 
$\rho \in I$, and since $\Phi''$ is $0$--homogeneous, recall \eqref{eq:Phidd}, 
we have from \eqref{eq:approx} and \eqref{eq:h1err} that
\begin{align*}
| S_1 | & \leq \max_{| p | =1} | \Phi''(p) | | x^m - x^m_h |_1 
\, | x^{m+1}_h - x^m_h |_1 \\ & 
\leq C \bigl( | x^m - Q_h x^m |_1 
+ | Q_h x^m - x^m_h |_1 \bigr) | x^{m+1}_h- x^m_h |_1  \\ & 
\leq C h | x^{m+1}_h - x^m_h |_1 
\leq \tfrac12 \hat{c} | x^{m+1}_h - x^m_h |_1^2 + C h^2.
\end{align*}
Clearly, 
\begin{displaymath}
|S_2 | \leq C \| x^{m+1}_h - x^m_h \|_0 \leq \frac{\sigma}{\Delta t} \| x^{m+1}_h - x^m_h \|_0^2 + C \Delta t \leq  \frac{\sigma}{\Delta t} \| x^{m+1}_h - x^m_h \|_0^2 + C \delta h,
\end{displaymath}
so that \eqref{eq:timedif} implies that
\begin{equation*} 
| x^{m+1}_h - x^m_h |_1 \leq C \bigl( h+ (\delta h)^{\frac12} \bigr).
\end{equation*}
The bound \eqref{eq:timedif1} then follows with the help of the inverse estimate 
\eqref{eq:inverse}.
\end{proof}

\begin{remark} \label{rem:abshom}
We remark that absolute 1-homogeneity assumption in \eqref{eq:gamma1}
is only used 
for the proof of Lemma~\ref{lem:ex} via the estimate \eqref{eq:Phiconvex3}.
We expect that the proof, and hence the results in this paper, can be extended
to positively homogeneous anisotropies, i.e.\
$\phi(\lambda p) = \lambda  \phi(p)$ for $p \in \bR^d$, 
$\lambda \in \bRplus$,
under a more restrictive time step condition.
\end{remark}

For later use we remark that 
in $I$ it holds that
\begin{equation} \label{eq:xmp1}
[x^m_{h,\rho},x^{m+1}_{h,\rho}], \, [x^{m+1}_\rho, x^{m+1}_{h,\rho}], \,
[(Q_h x^m)_\rho,(Q_h x^{m+1})_\rho], \,
[(Q_h x^m)_\rho, x^{m+1}_{h,\rho}], \, [(Q_h x^{m+1})_\rho, x^{m+1}_{h,\rho}]
 \subset K, 
\end{equation}
provided that $h_0$ and $\delta$ are sufficiently small.
For example, if $z=(1-\lambda) x^m_{h,\rho} + \lambda x^{m+1}_{h,\rho} \in [x^m_{h,\rho},x^{m+1}_{h,\rho}]$, then we have similarly to \eqref{eq:z} that
\begin{displaymath}
| z | \geq | x^m_\rho | - | x^m_\rho - x^m_{h,\rho} | - \lambda | x^m_{h,\rho} - x^{m+1}_{h,\rho} | \geq c_0 - \tfrac12 c_0 - C(h^{\frac{1}{2}}+ \delta^{\frac{1}{2}})
\geq \tfrac14 c_0,
\end{displaymath}
where we have used \eqref{eq:regul}, \eqref{eq:ind} and \eqref{eq:timedif1}.
The other inclusions can be shown in a similar way, on also making use of
\eqref{eq:xregul}.
In particular we obtain in a similar way as in \eqref{eq:h1err} that
\begin{equation} \label{eq:h1err1}
| e^{m+1}_h |_1^2  \leq C \bigl( E^{m+1}+F^{m+1}  \bigr) + C h^4.
\end{equation}

Evaluating \eqref{eq:Hxtweak} for $t=t_{m+1}$ we have
\begin{displaymath}
\int_I H(x^{m+1}_\rho) x^{m+1}_t \cdot \eta \drho 
+ \int_I \Phi'(x^{m+1}_\rho) \cdot \eta_\rho \drho =0 
\qquad \forall\ \eta \in [H^1(I)]^d.
\end{displaymath}
Combining this relation with \eqref{eq:fea}, 
and recalling \eqref{eq:defeh}, \eqref{eq:defxizeta} and \eqref{eq:defhatx}, 
we obtain the error equation
\begin{align*} 
& \frac1{\Delta t} \int_I H(x^m_{h,\rho}) (e^{m+1}_h-e^m_h) \cdot \eta_h \drho
 + \int_I \bigl( \Phi'((Q_h x^{m+1})_\rho)
- \Phi'(x^{m+1}_{h,\rho}) \bigr) \cdot \eta_{h,\rho} \drho \nonumber \\ & \quad
= \int_I \bigl( H(x^m_{h,\rho}) - H(x^{m+1}_\rho) \bigr) x^{m+1}_t 
\cdot \eta_h \drho
+ \int_I H(x^m_{h,\rho}) 
\left[ \frac{Q_h x^{m+1}- Q_h x^m}{\Delta t} - x^{m+1}_t \right] \cdot \eta_h 
\drho \nonumber \\ & \qquad 
+ \int_I \zeta^{m+1} \cdot \eta_{h,\rho} 
+ \int_I \Phi''(x^{m+1}_\rho) ( Q_h x^{m+1} - x^{m+1})_\rho \cdot \eta_{h,\rho}
\drho \nonumber \\ & \quad
= \int_I \bigl( H(x^m_{h,\rho}) - H(x^{m+1}_\rho) \bigr) x^{m+1}_t 
\cdot \eta_h \drho
+ \int_I H(x^m_{h,\rho}) 
\left[ \frac{Q_h x^{m+1}- Q_h x^m}{\Delta t} - x^{m+1}_t \right] \cdot \eta_h 
\drho \nonumber \\ & \qquad 
+ \int_I \zeta^{m+1} \cdot \eta_{h,\rho} 
+ \int_I  ( x^{m+1}-Q_h x^{m+1})\cdot \eta_h \drho.
\end{align*}
Choosing $\eta_h=e^{m+1}_h-e^m_h$ and recalling \eqref{eq:Hposdef1} we obtain
\begin{align} \label{eq:err2} 
& \frac{\sigma}{\Delta t} \| e^{m+1}_h - e^m_h \|_0^2+ \int_I \bigl( \Phi'((Q_h x^{m+1})_\rho)
- \Phi'(x^{m+1}_{h,\rho}) \bigr) \cdot (e^{m+1}_{h,\rho} - e^m_{h,\rho}) \drho 
 \nonumber \\ & \quad
 \leq  \int_I \bigl( H(x^m_{h,\rho}) - H(x^{m+1}_\rho) \bigr) x^{m+1}_t 
\cdot ( e^{m+1}_h - e^m_h) \drho \nonumber \\ & \qquad 
+ \int_I H(x^m_{h,\rho}) \left[ \frac{Q_h x^{m+1}- Q_h x^m}{\Delta t} - x^{m+1}_t \right] \cdot (e^{m+1}_h - e^m_h) \drho 
\nonumber \\ & \qquad
 + \int_I \zeta^{m+1} \cdot (e^{m+1}_{h,\rho} - e^m_{h,\rho})  \drho  
+ \int_I  ( x^{m+1}-Q_h x^{m+1})\cdot (e^{m+1}_h - e^m_h) \drho
\nonumber \\ & \quad 
=: T_1+T_2 + T_3+T_4.
\end{align}
The treatment of the second term on the left hand side of \eqref{eq:err2} 
is quite involved. That is why we deal with it in the following lemma.

\newcommand{\Qhx}{{\mathfrak X}}
\begin{lemma} \label{lem:lhs}
It holds that
\begin{align*}
& \lefteqn{ \int_I \bigl( \Phi'((Q_h x^{m+1})_\rho)
- \Phi'(x^{m+1}_{h,\rho}) \bigr) \cdot (e^{m+1}_{h,\rho} - e^m_{h,\rho}) \drho } \\ & \quad 
\geq E^{m+1} - E^m + \tfrac14 \sigma | e^{m+1}_h -e^m_h |_1^2 - C (\Delta t)^3 - C \Delta t  \bigl( | e^m_h |_1^2 + | e^{m+1}_h |_1^2 \bigr),
\end{align*}
provided that $h_0$ and $\delta$ are small enough.
\end{lemma}
\begin{proof} 
Throughout this proof we use the shorthand notations $\Qhx$ for
$Q_h x$ and $\Qhx^m$ for $Q_h x^m$.
Let
\begin{align*}
D^m &= \bigl( \Phi'(\Qhx^{m+1}_\rho) - \Phi'(x^{m+1}_{h,\rho}) \bigr) \cdot (e^{m+1}_{h,\rho} - e^m_{h,\rho}) \\
& = \bigl( \Phi'(\Qhx^{m+1}_\rho) - \Phi'(x^{m+1}_{h,\rho}) \bigr) \cdot \bigl( (\Qhx^{m+1}_\rho - \Qhx^m_\rho) - (x^{m+1}_{h,\rho} 
- x^m_{h,\rho}) \bigr)  \\
& = \bigl( \Phi(x^{m+1}_{h,\rho}) - \Phi'(\Qhx^{m+1}_\rho) \cdot x^{m+1}_{h,\rho} \bigr) - 
\bigl( \Phi(x^m_{h,\rho}) - \Phi'(\Qhx^m_\rho) \cdot x^m_{h,\rho} \bigr) 
\\ & \,
+ [\Phi(x^m_{h,\rho}) - \Phi(x^{m+1}_{h,\rho}) - \Phi'(x^{m+1}_{h,\rho}) \cdot (x^m_{h,\rho}- x^{m+1}_{h,\rho})]
+ \bigl( \Phi'(\Qhx^{m+1}_\rho) - \Phi'(\Qhx^m_\rho) \bigr) \cdot x^m_{h,\rho} \\ & \,
+ \Phi'(\Qhx^{m+1}_\rho) \cdot (\Qhx^{m+1}_\rho - \Qhx^m_\rho) - \Phi'(x^{m+1}_{h,\rho}) \cdot (\Qhx^{m+1}_\rho -\Qhx^m_\rho).
\end{align*}
Since $\Phi'(p) \cdot p = 2 \Phi(p)$, we can rewrite this as
\begin{align}
D^m & = \bigl[ \Phi(x^{m+1}_{h,\rho}) - \Phi(\Qhx^{m+1}_\rho) -  \Phi'(\Qhx^{m+1}_\rho) \cdot (x^{m+1}_{h,\rho} -\Qhx^{m+1}_\rho) \bigr]
\nonumber \\ & \quad 
- \bigl[ \Phi(x^m_{h,\rho}) - \Phi(\Qhx^m_\rho) -  \Phi'(\Qhx^m_\rho) \cdot (x^m_{h,\rho}-\Qhx^m_\rho) \bigr] + D^m_1+D^m_2+D^m_3, \label{eq:err3}
\end{align}
where
\begin{align*}
D^m_1 &= \Phi(\Qhx^m_\rho) - \Phi(\Qhx^{m+1}_\rho) - \Phi'(\Qhx^{m+1}_\rho) \cdot (\Qhx^m_\rho- \Qhx^{m+1}_\rho), \\
D^m_2 &=  \Phi(x^m_{h,\rho}) - \Phi(x^{m+1}_{h,\rho}) - \Phi'(x^{m+1}_{h,\rho}) \cdot (x^m_{h,\rho}- x^{m+1}_{h,\rho}), \\
D^m_3 &=   \bigl( \Phi'(\Qhx^{m+1}_\rho) - \Phi'(\Qhx^m_\rho) \bigr) \cdot x^m_{h,\rho} 
- \Phi'(x^{m+1}_{h,\rho}) \cdot (\Qhx^{m+1}_\rho -\Qhx^m_\rho).
\end{align*}
Using a Taylor expansion, we obtain 
\begin{align*}
D^m_1 &= \tfrac12 \Phi''(\Qhx^{m+1}_\rho)(\Qhx^m_\rho - \Qhx^{m+1}_\rho) \cdot (\Qhx^m_\rho - \Qhx^{m+1}_\rho) \\ & \quad
+ \tfrac16 \Phi'''(\theta_1)(\Qhx^m_\rho - \Qhx^{m+1}_\rho, \Qhx^m_\rho - \Qhx^{m+1}_\rho, \Qhx^m_\rho - \Qhx^{m+1}_\rho) \\ & \
= \tfrac12 \Phi''(x^{m+1}_{h,\rho})(\Qhx^m_\rho - \Qhx^{m+1}_\rho) \cdot (\Qhx^m_\rho - \Qhx^{m+1}_\rho) \\ & \quad
+ \tfrac12 [\Phi''(\Qhx^{m+1}_\rho) - \Phi''(x^{m+1}_{h,\rho})] 
(\Qhx^m_\rho - \Qhx^{m+1}_\rho) \cdot (\Qhx^m_\rho - \Qhx^{m+1}_\rho) \\ & \quad
+ \tfrac16 \Phi'''(\theta_1)(\Qhx^m_\rho - \Qhx^{m+1}_\rho, \Qhx^m_\rho - \Qhx^{m+1}_\rho, \Qhx^m_\rho - \Qhx^{m+1}_\rho) \\ & \
=: \tfrac12 \Phi''(x^{m+1}_{h,\rho})(\Qhx^m_\rho - \Qhx^{m+1}_\rho) \cdot (\Qhx^m_\rho - \Qhx^{m+1}_\rho) + R^m_1,
\end{align*}
for some $\theta_1 \in [\Qhx^m_\rho,\Qhx^{m+1}_\rho]
\subset K$, recall \eqref{eq:xmp1}. 
Since $| \Qhx_t |_{1,\infty} \leq C$ by  \eqref{eq:approxtime} and \eqref{eq:xregul},
we infer on recalling \eqref{eq:xmp1} and \eqref{eq:defeh} that
\begin{align} \label{eq:r1}
| R^m_1 | & \leq C | \Qhx^{m+1}_\rho - \Qhx^m_\rho |^2 \bigl( | \Qhx^{m+1}_\rho - \Qhx^m_\rho | + | \Qhx^{m+1}_\rho - x^{m+1}_{h,\rho}| \bigr)
\nonumber \\ &
\leq C (\Delta t)^2 \bigl( \Delta t + | e^{m+1}_{h,\rho} | \bigr).
\end{align}
Similarly, we obtain 
\begin{align*}
D^m_2 = \tfrac12 \Phi''(x^{m+1}_{h,\rho})(x^m_{h,\rho} - x^{m+1}_{h,\rho}) \cdot (x^m_{h,\rho} - x^{m+1}_{h,\rho}) + R^m_2,
\end{align*}
where, recalling \eqref{eq:xmp1}, 
\begin{equation} \label{eq:r2}
| R^m_2 | \leq C | x^m_{h,\rho} - x^{m+1}_{h,\rho} |^3.
\end{equation}
If we combine the expressions for $D^m_1$ and $D^m_2$ with $D^m_3$
and recall \eqref{eq:defeh}, we obtain
\begin{align}
D^m_1+D^m_2+D^m_3 & = 
 \tfrac12 \Phi''(x^{m+1}_{h,\rho}) (e^{m+1}_{h,\rho} - e^m_{h,\rho}) \cdot (e^{m+1}_{h,\rho} - e^m_{h,\rho}) \nonumber \\ & \quad 
+ \Phi''(x^{m+1}_{h,\rho}) (x^{m+1}_{h,\rho}-x^m_{h,\rho}) \cdot (\Qhx^{m+1}_\rho - \Qhx^m_\rho) + R^m_1 + R^m_2  \nonumber \\ & \quad
+ \bigl( \Phi'(\Qhx^{m+1}_\rho) - \Phi'(\Qhx^m_\rho) \bigr) \cdot x^m_{h,\rho} 
- \Phi'(x^{m+1}_{h,\rho}) \cdot (\Qhx^{m+1}_\rho -\Qhx^m_\rho) \nonumber \\ & 
=  \tfrac12 \Phi''(x^{m+1}_{h,\rho}) (e^{m+1}_{h,\rho} - e^m_{h,\rho}) \cdot (e^{m+1}_{h,\rho} - e^m_{h,\rho}) +  R^m_1+ R^m_2 + R^m_3, \label{eq:err4}
\end{align}
where we used \eqref{eq:Phidd} in the last step and defined
\begin{equation*} 
R^m_3 = \bigl( \Phi'(\Qhx^{m+1}_\rho) - \Phi'(\Qhx^m_\rho) - \Phi''(x^{m+1}_{h,\rho}) (\Qhx^{m+1}_\rho - \Qhx^m_\rho) \bigr) \cdot x^m_{h,\rho}.
\end{equation*}
In order to estimate $R^m_3$ we again use a Taylor expansion and obtain
\begin{align*}
R^m_3 & = \bigl( \Phi''(\Qhx^m_\rho) - \Phi''(x^{m+1}_{h,\rho})\bigr) 
(\Qhx^{m+1}_\rho - \Qhx^m_\rho)  \cdot x^m_{h,\rho} \nonumber \\ & \quad
+ \tfrac12 \Phi'''(\theta_2)(\Qhx^{m+1}_\rho - \Qhx^m_\rho,\Qhx^{m+1}_\rho - \Qhx^m_\rho,x^m_{h,\rho}) \\ & 
= \Phi'''(\theta_3)(\Qhx^m_\rho- x^{m+1}_{h,\rho}, 
\Qhx^{m+1}_\rho - \Qhx^m_\rho,x^m_{h,\rho}) \\ & \quad
+ \tfrac12 \Phi'''(\theta_2)(\Qhx^{m+1}_\rho - \Qhx^m_\rho,\Qhx^{m+1}_\rho - \Qhx^m_\rho,x^m_{h,\rho}),
\end{align*}
where $\theta_2 \in [\Qhx^m_\rho,\Qhx^{m+1}_\rho] \subset K$, 
$\theta_3 \in [\Qhx^m_\rho,x^{m+1}_{h,\rho}]\subset K$, 
recall again \eqref{eq:xmp1}.  
Then, on recalling \eqref{eq:Phidd} and \eqref{eq:approxtime}, we have that
\begin{align} \label{eq:r3}
| R^m_3 | & \leq 
| \Phi'''(\theta_3)(\Qhx^m_\rho- x^{m+1}_{h,\rho},\Qhx^{m+1}_\rho - \Qhx^m_\rho,x^m_{h,\rho}- \theta_3) | \nonumber \\ & \qquad
 + \tfrac12 | \Phi'''(\theta_2)(\Qhx^{m+1}_\rho - \Qhx^m_\rho,\Qhx^{m+1}_\rho - \Qhx^m_\rho,x^m_{h,\rho}-\theta_2) | \nonumber \\ & 
\leq C \Delta t | \Qhx^m_\rho - x^{m+1}_{h,\rho} | | x^m_{h,\rho} - \theta_3| + C (\Delta t)^2 | x^m_{h,\rho} - \theta_2 | \nonumber \\ & 
\leq C \Delta t \bigl( \Delta t + | e^{m+1}_{h,\rho} | \bigr) \bigl( | e^m_{h,\rho}| +
| x^{m+1}_{h,\rho} - x^m_{h,\rho} | \bigr) + C (\Delta t)^2 \bigl( \Delta t + | e^m_{h,\rho} | \bigr) \nonumber \\ & 
\leq C (\Delta t)^3 + C \Delta t \bigl( | e^m_{h,\rho} |^2 + | e^{m+1}_{h,\rho} |^2 \bigr) +  C | x^{m+1}_{h,\rho} - x^m_{h,\rho} |^3,
\end{align}
where we have observed that for a $\lambda \in [0,1]$ we can write
\begin{align*}
|x^m_{h,\rho} - \theta_3| &
= |x^m_{h,\rho} - \lambda \Qhx^m_\rho - (1-\lambda) x^{m+1}_{h,\rho}|
\leq \lambda | x^m_{h,\rho} - \Qhx^m_\rho | + (1-\lambda) 
| x^m_{h,\rho} - x^{m+1}_{h,\rho}| \\ &
\leq | e^m_{h,\rho}| + | x^m_{h,\rho} - x^{m+1}_{h,\rho}|.
\end{align*}
Combining \eqref{eq:err4}, \eqref{eq:r1}, \eqref{eq:r2} and \eqref{eq:r3}, 
and recalling \eqref{eq:convex3} and \eqref{eq:xmp1}, then yields
\begin{align*} 
D^m_1+D^m_2+D^m_3 & \geq  \tfrac12 \sigma | e^{m+1}_{h,\rho} - e^m_{h,\rho} |^2 
- C (\Delta t)^3 
- C \Delta t \bigl( | e^m_{h,\rho} |^2 + | e^{m+1}_{h,\rho} |^2 \bigr) 
- C | x^{m+1}_{h,\rho} - x^{m}_{h,\rho} |^3 .
\end{align*}
Inserting this estimate into \eqref{eq:err3},
integrating the resulting inequality with respect to $\rho$ and recalling the definitions of $D^m$ and $E^m$ yields
\begin{align}
& \int_I \bigl( \Phi'(\Qhx^{m+1}_\rho)
- \Phi'(x^{m+1}_{h,\rho}) \bigr) \cdot (e^{m+1}_{h,\rho} - e^m_{h,\rho}) \drho 
\nonumber \\ &
\geq E^{m+1} - E^m + \tfrac12 \sigma | e^{m+1}_h -e^m_h |_1^2 - C (\Delta t)^3 - C \Delta t  \bigl( | e^m_h |_1^2 + | e^{m+1}_h |_1^2 \bigr) - C \int_I | x^{m+1}_{h,\rho} - x^m_{h,\rho} |^3 \drho.  \label{eq:err5}
\end{align}
On recalling \eqref{eq:timedif1}, we have
\begin{align*} 
\int_I | x^{m+1}_{h,\rho} - x^m_{h,\rho} |^3 \drho & 
\leq 2 \int_I | x^{m+1}_{h,\rho} - x^m_{h,\rho} |
 \bigl( | \Qhx^{m+1}_\rho - \Qhx^m_\rho |^2 + | e^{m+1}_{h,\rho} - e^m_{h,\rho} |^2 \bigr) \drho  \\ & 
\leq C (\Delta t)^2 |x^{m+1}_h - x^m_h|_1 + C \bigl( h^{\frac12} + \delta^{\frac12} \bigr)  |e^{m+1}_h - e^m_h|_1^2 
\\ & 
\leq C (\Delta t)^2 \bigl( |  e^m_h |_1 + | e^{m+1}_h |_1 + \Delta t \bigr) + 
C \bigl( h^{\frac12} + \delta^{\frac12} \bigr)  | e^{m+1}_h - e^m_h |_1^2  \\ & 
\leq C (\Delta t)^3 + C \Delta t \bigl( | e^m_h |_1^2 + | e^{m+1}_h |_1^2\bigr)
 + \tfrac14 \sigma  | e^{m+1}_h - e^m_h |_1^2
\end{align*}
provided that $h_0$ and $\delta$ are small enough. Inserting the above relation into \eqref{eq:err5} finishes the proof.
\end{proof} 

Let us next turn to the terms on the right hand side of \eqref{eq:err2}. We have
\begin{lemma} \label{lem:rhs}
For every $\epsilon>0$ there exists a $C_\epsilon>0$ such that
\[
 \sum_{k=1}^4 T_k  \leq -F^{m+1} + F^m + \frac{\epsilon}{\Delta t} \| e^{m+1}_h - e^m_h \|_0^2 + C_\epsilon \Delta t  \bigl( h^4 + (\Delta t)^2 \bigr) 
+ C_\epsilon \Delta t \bigl( | e^m_h |_1^2 + | e^{m+1}_h |_1^2 \bigr). 
\]
\end{lemma}
\begin{proof}
We again use the shorthand notations $\Qhx$ for 
$Q_h x$ and $\Qhx^m$ for $Q_h x^m$.
To begin, we write
\begin{align*}
T_1 & = \int_I \bigl( H(x^m_{h,\rho}) - H(\Qhx^m_\rho) \bigr) x^{m+1}_t  \cdot ( e^{m+1}_h - e^m_h) \drho \nonumber \\ & \quad 
+ \int_I \bigl( H(\Qhx^m_\rho) - H(x^m_\rho) \bigr) x^{m+1}_t  \cdot ( e^{m+1}_h - e^m_h)  \drho 
\nonumber \\ & \quad 
+ \int_I \bigl( H(x^m_\rho) - H(x^{m+1}_\rho) \bigr) x^{m+1}_t  \cdot ( e^{m+1}_h - e^m_h) \drho 
=: T_{1,1} +  T_{1,2} + T_{1,3}. 
\end{align*}
Clearly, it holds that
\begin{align}
| T_{1,1} | + | T_{1,3} | & \leq 
C \bigl( | e^m_h |_1 + | x^m - x^{m+1} |_1 \bigr) \| e^{m+1}_h - e^m_h \|_0
\nonumber  \\ & 
\leq \frac{\epsilon}{\Delta t} \| e^{m+1}_h - e^m_h \|_0^2 
+ C_\epsilon \Delta t | e^m_h |_1^2 + C_\epsilon (\Delta t)^3.
\label{eq:t1113}
\end{align}
In order to treat $T_{1,2}$ we write, on recalling \eqref{eq:defxizeta},
\[
(H(\Qhx^m_\rho) - H(x^m_\rho)) x^{m+1}_t = 
( \Qhx^m_{i,\rho} - x^m_{i,\rho})\xi^m_i  + R^m_4,
\  \text{ where } \ | R^m_4 | \leq C | \Qhx^m_\rho - x^m_\rho |^2.
\]
Then, with the help of \eqref{eq:approx}, \eqref{eq:approxtime} and
\eqref{eq:xregul}, we have that
\begin{align}
T_{1,2} & = \int_I ( \Qhx^m_{i,\rho} - x^m_{i,\rho}) \xi^m_i \cdot (e^{m+1}_h - e^m_h) \drho + \int_I R^m_4 \cdot
(e^{m+1}_h - e^m_h) \drho \nonumber \\ &
= - \int_I ( \Qhx^m_i - x^m_i) \xi^m_{i,\rho} \cdot (e^{m+1}_h - e^m_h ) \drho - \int_I 
 ( \Qhx^m_i - x^m_i) \xi^m_i \cdot (e^{m+1}_{h,\rho} - e^m_{h,\rho} ) \drho 
\nonumber \\ & \quad
+ \int_I R^m_4 \cdot (e^{m+1}_h - e^m_h) \drho \nonumber \\ & 
= - \int_I ( \Qhx^{m+1}_i - x^{m+1}_i) \xi^{m+1}_i \cdot e^{m+1}_{h,\rho} \drho +
 \int_I ( \Qhx^m_i - x^m_i) \xi^m_i \cdot e^m_{h,\rho} \drho 
\nonumber \\ & \quad
+ \int_I \left[  \bigl( (\Qhx^{m+1}_i - \Qhx^m_i) - (x^{m+1}_i - x^m_i) \bigr) \xi^m_i \right]  \cdot e^{m+1}_{h,\rho} \drho 
\nonumber \\ & \quad
+ \int_I \left[
 (\Qhx^{m+1}_i - x^{m+1}_i) ( \xi^{m+1}_i - \xi^m_i) \right]  \cdot e^{m+1}_{h,\rho} \drho 
\nonumber \\ & \quad
- \int_I ( \Qhx^m_i - x^m_i) \xi^m_{i,\rho} \cdot (e^{m+1}_h - e^m_h ) \drho + \int_I R^m_4 \cdot (e^{m+1}_h - e^m_h) \drho 
\nonumber \\ & 
\leq - \int_I ( \Qhx^{m+1}_i - x^{m+1}_i) \xi^{m+1}_i \cdot e^{m+1}_{h,\rho} \drho +
 \int_I ( \Qhx^m_i - x^m_i) \xi^m_i \cdot e^m_{h,\rho} \drho 
\nonumber \\ & \quad
+ C  \Delta t h^2 | e^{m+1}_h |_1 + C h^2 \| e^{m+1}_h - e^m_h \|_0 
\nonumber \\ & 
\leq - \int_I ( \Qhx^{m+1}_i - x^{m+1}_i) \xi^{m+1}_i \cdot e^{m+1}_{h,\rho} \drho +
 \int_I ( \Qhx^m_i - x^m_i) \xi^m_i \cdot e^m_{h,\rho} \drho 
\nonumber \\ & \quad
+ \frac{\epsilon}{\Delta t} \| e^{m+1}_h - e^m_h \|_0^2 
+ C_\epsilon \Delta t  \bigl( h^4 + | e^{m+1}_h |_1^2 \bigr). \label{eq:t12}
\end{align}
Combining \eqref{eq:t1113} and \eqref{eq:t12} we have
\begin{align}
T_1 & \leq - \int_I ( \Qhx^{m+1}_i - x^{m+1}_i) \xi^{m+1}_i \cdot e^{m+1}_{h,\rho} \drho +
 \int_I ( \Qhx^m_i - x^m_i) \xi^m_i \cdot e^m_{h,\rho} \drho 
\nonumber \\ & \quad
+ \frac{2\epsilon}{\Delta t} \| e^{m+1}_h - e^m_h \|_0^2 
+ C_\epsilon \Delta t  \bigl( h^4 + (\Delta t)^2 \bigr) 
+ C_\epsilon \Delta t \bigl( | e^m_h |_1^2 + | e^{m+1}_h |_1^2 \bigr). 
\label{eq:t1}
\end{align}
Next, on recalling \eqref{eq:approxtime}, we obtain 
\begin{align}
| T_2 | & 
\leq C \left\| \frac{\Qhx^{m+1} - \Qhx^m}{\Delta t} - x^{m+1}_t \right\|_0 
\| e^{m+1}_h - e^m_h \|_0 \nonumber  \\ & 
= C \frac{1}{\Delta t} \left\| \int_{t_m}^{t_{m+1}} (\Qhx_t - x^{m+1}_t) \dt 
\right\|_0  \| e^{m+1}_h - e^m_h \|_0 \nonumber  \\ & 
\leq C \bigl( \frac{1}{\Delta t} \int_{t_m}^{t_{m+1}} \| \Qhx_t - x_t \|_0
 \dt + \Delta t  \sup_{0 < t < T} \Vert x_{tt} \Vert_0 \bigr) \| e^{m+1}_h - e^m_h \|_0 \nonumber \\ &
\leq C \bigl( h^2 + \Delta t \bigr) \| e^{m+1}_h - e^m_h \|_0 
\leq \frac{\epsilon}{\Delta t} \| e^{m+1}_h - e^m_h \|_0^2 
+ C_\epsilon \Delta t \bigl( h^4 + (\Delta t)^2 \bigr). \label{eq:t2}
\end{align}
For the third term on the right hand side of \eqref{eq:err2} we have
\begin{align}
T_3 & 
 = \int_I \zeta^{m+1} \cdot e^{m+1}_{h,\rho} \drho - \int_I \zeta^m \cdot e^m_{h,\rho} \drho -
 \int_I ( \zeta^{m+1}-\zeta^m) \cdot e^m_{h,\rho} \drho \nonumber \\ &
= \int_I \zeta^{m+1} \cdot e^{m+1}_{h,\rho} \drho - \int_I \zeta^m \cdot e^m_{h,\rho} \drho -
 \int_{t^m}^{t_{m+1}} \int_I \zeta_t \cdot e^m_{h,\rho} \drho \dt. \label{eq:t3a}
\end{align}
We have from \eqref{eq:defxizeta} and \eqref{eq:Phidd} that 
$\zeta = \Phi'(\Qhx_\rho) - \Phi''(x_\rho) \Qhx_\rho$, and so
\begin{align*}
\zeta_t & = (\Phi''(\Qhx_\rho) - \Phi''(x_\rho))  \Qhx_{\rho t} - \Phi'''(x_\rho) (\Qhx_\rho- x_\rho,x_{\rho t},\cdot) 
\\ & 
= \int_0^1 \Phi'''(s \Qhx_\rho+(1-s) x_\rho) (\Qhx_\rho - x_\rho, \Qhx_{\rho t},\cdot) \ds
 - \Phi'''(x_\rho) (\Qhx_\rho- x_\rho,x_{\rho t},\cdot) \\ & 
= \int_0^1 \left[ \Phi'''(s \Qhx_\rho+(1-s) x_\rho)- \Phi'''(x_\rho)\right]
 (\Qhx_\rho - x_\rho, \Qhx_{\rho t},\cdot) \ds \\ & \quad
 + \Phi'''(x_\rho)(\Qhx_\rho-x_\rho, \Qhx_{\rho t} - x_{\rho t},\cdot),
\end{align*}
so that 
\begin{equation} \label{eq:zetatest}
| \zeta_t | \leq C | \Qhx_\rho - x_\rho | 
\bigl(  | \Qhx_\rho - x_\rho | + | \Qhx_{\rho t} - x_{\rho t} | \bigr).
\end{equation}
Combining \eqref{eq:t3a} and \eqref{eq:zetatest} yields, on recalling
\eqref{eq:approxboth} and \eqref{eq:approxtime}, that 
\begin{align} \label{eq:t3b}
T_3 & \leq \int_I \zeta^{m+1} \cdot e^{m+1}_{h,\rho} \drho - \int_I \zeta^m \cdot e^m_{h,\rho} \drho \nonumber \\ & \quad
+ C \Delta t \sup_{t_m \leq t \leq t_{m+1}} 
| \Qhx - x |_{1,\infty} 
\sup_{t_m \leq t \leq t_{m+1}} 
\bigl( | \Qhx - x |_1 + | \Qhx_t - x_t |_1 \bigr) | e^m_h |_1 \nonumber\\ &
\leq \int_I \zeta^{m+1} \cdot e^{m+1}_{h,\rho} \drho 
- \int_I \zeta^m \cdot e^m_{h,\rho} \drho + C \Delta t h^4
+ C \Delta t | e^m_h |_1^2. 
\end{align}
Lastly, we can bound
\begin{equation} \label{eq:t4}
| T_4 | \leq  \| x^{m+1} - \Qhx^{m+1} \|_0 \| e^{m+1}_h - e^m_h \|_0
\leq C h^2 \| e^{m+1}_h - e^m_h \|_0 
\leq \frac{\epsilon}{\Delta t} \| e^{m+1}_h - e^m_h \|_0^2 
+ C_\epsilon \Delta t h^4.
\end{equation}
Collecting   \eqref{eq:t1}, \eqref{eq:t2}, \eqref{eq:t3b} and \eqref{eq:t4} and recalling the definition of $F^m$ yields the result.
\end{proof}

Let us now insert the estimates obtained in  Lemma \ref{lem:lhs} and Lemma \ref{lem:rhs} into \eqref{eq:err2}.  After choosing 
$\epsilon$ sufficiently small we obtain
\begin{align} \label{eq:err7}
& 
\frac{\sigma}{2 \Delta t} \| e^{m+1}_h - e^m_h \|_0^2 + (E^{m+1} + F^{m+1}) - (E^m+ F^m) 
+ \tfrac14 \sigma | e^{m+1}_h - e^m_h |_1^2 
\nonumber \\ & \quad
\leq C \Delta t \bigl( h^4 + (\Delta t)^2 \bigr) 
+ C \Delta t \bigl( | e^m_h |_1^2 + | e^{m+1}_h |_1^2 \bigr) 
\nonumber \\ & \quad
\leq C_3 \Delta t \bigl( h^4 + (\Delta t)^2 \bigr) 
+ C_3 \Delta t \bigl( E^m+ F^m + E^{m+1}+ F^{m+1} \bigr) ,
\end{align}
in view of \eqref{eq:h1err} and \eqref{eq:h1err1}. 
If we choose $h_0$ so small that 
$C_3 \Delta t \leq \frac12$ for $\Delta t \leq \delta h_0$, then
$0 < (1 - C_3 \Delta t)^{-1} \leq 1 + 2 C_3 \Delta t$, and so it follows from
\eqref{eq:ind} that
\begin{align*}
E^{m+1}+F^{m+1} & \leq (1 - C_3 \Delta t)^{-1} \left[(1+ C_3 \Delta t)(E^m+F^m) 
+ C_3 \Delta t \bigl( h^4 + (\Delta t)^2 \bigr) \right] \\ &
\leq \bigl( 1 + 2 C_3 \Delta t \bigr)^2 (E^m+F^m) + C_3 \bigl( 1 + 2C_3 \Delta t \bigr) \Delta t \bigl( h^4 + (\Delta t)^2 \bigr) 
\\ & 
\leq \bigl( 1 + 2 C_3 \Delta t \bigr)^2 \bigl( h^4 + (\Delta t)^2 \bigr) e^{\mu t_m} + 2C_3\Delta t \bigl( h^4 + (\Delta t)^2 \bigr) 
\\ & 
\leq \bigl( 1 + 3 C_3 \Delta t \bigr)^2 \bigl( h^4 + (\Delta t)^2 \bigr) e^{\mu t_m} \\ &
\leq \bigl( h^4 + (\Delta t)^2 \bigr) e^{6C_3 \Delta t} e^{\mu t_m} = \bigl( h^4 + (\Delta t)^2 \bigr) e^{\mu t_{m+1}},
\end{align*}
if we choose $\mu = 6C_3$. This proves the second estimate in \eqref{eq:ind}. 
In order to show the first estimate in \eqref{eq:ind}, we observe from
\eqref{eq:defeh}, \eqref{eq:inverse}, \eqref{eq:approxinf}
and \eqref{eq:h1err} with $m$ replaced by $m+1$, that
\begin{displaymath}
| x^{m+1}_{h,\rho} - x^{m+1}_\rho| 
\leq | e^{m+1}_h |_{1,\infty} + | Q_h x^{m+1} - x^{m+1} |_{1,\infty}
\leq C h^{-\frac12} | e^{m+1}_h |_1 + C h \leq C h^\frac12
\leq \tfrac12 c_0,
\end{displaymath}
provided that $h_0$ is chosen sufficiently small.
Since $\mu$ and $\delta$, were chosen independently of $h$ and 
$\Delta t$, we have shown \eqref{eq:ind} by induction. 

It remains to show that \eqref{eq:ind} implies the desired result 
\eqref{eq:eb}. 
The second bound in \eqref{eq:eb} follows from \eqref{eq:approx}
and \eqref{eq:h1err}, since $x^m- x^m_h=x^m - Q_h x^m+ e^m_h$. 
In order to prove the $L^2$--error bound in \eqref{eq:eb}, 
we first remark that \eqref{eq:err7} together with \eqref{eq:ind} implies that
\begin{equation} \label{eq:err8}
\frac{1}{\Delta t} \sum_{m=0}^{M-1} \Vert e^{m+1}_h - e^m_h \Vert_0^2 \leq C \bigl( h^4+ (\Delta t)^2 \bigr).
\end{equation}
Since $e^0_h=0$, we obtain with the help of \eqref{eq:err8} for $1 \leq \ell \leq M$
\begin{align*}
\Vert e^\ell_h \Vert_0^2 & = \sum_{m=0}^{\ell-1} \bigl( \Vert e^{m+1}_h \Vert_0^2 - \Vert e^m_h \Vert_0^2 \bigr) 
\leq \sum_{m=0}^{\ell-1} \left(\Vert e^{m+1}_h - e^m_h \Vert_0 \bigl( \Vert e^{m+1}_h \Vert_0+ \Vert e^m_h \Vert_0 \bigr)\right) \\
& \leq 2 \Bigl( \frac{1}{\Delta t} \sum_{m=0}^{M-1} \Vert e^{m+1}_h - e^m_h \Vert_0^2 \Bigr)^{\frac{1}{2}} \Bigl( \Delta t \sum_{m=0}^{\ell}   
\Vert e^m_h \Vert_0^2  \Bigr)^{\frac{1}{2}} 
\leq C \bigl( h^4 + (\Delta t)^2 \bigr) + \Delta t \sum_{m=0}^{\ell} \Vert e^m_h \Vert_0^2.
\end{align*}
The discrete Gronwall lemma yields that $\max_{0 \leq m \leq M} \Vert e^m_h \Vert_0^2 \leq C(h^4 + (\Delta t)^2)$, so that the $L^2$--bound in \eqref{eq:eb} now follows
again from \eqref{eq:approx}. 

\setcounter{equation}{0}
\section{Numerical results} \label{sec:nr}

We implemented the scheme \eqref{eq:fea} within the 
finite element toolbox Alberta, \cite{Alberta}. The systems of nonlinear
equations arising at each time level are solved using a Newton
iteration, where the solutions to the linear subsystems are found
with the help of the sparse factorization package UMFPACK, see \cite{Davis04}.
For all our numerical simulations we use a uniform partitioning of $[0,1]$,
so that $q_j = jh$, $j=0,\ldots,J$, with $h = \frac 1J$. Unless otherwise
stated, we use $J=512$ and $\Delta t = 10^{-4}$.
For our numerical simulations we will often be interested in the evolution
of the discrete energy $E_\phi(x^m_h) 
= \int_I \phi(x^m_{h,\rho}) \drho$, recall \eqref{eq:Ephi}. 
We also consider the ratio
\begin{equation*} 
\ratio^m = \dfrac{\max_{j=1,\ldots,J} |x_h^m(q_j) - x_h^m(q_{j-1})|}
{\min_{j=1,\ldots, J} |x_h^m(q_j) - x_h^m(q_{j-1})|}
\end{equation*}
between the longest and shortest element of $\Gamma^m_h = x_h^m(I)$, 
and are often interested in the evolution of this ratio over time.
We stress that no redistribution of vertices was necessary during any of our
numerical simulations.

Moreover, we will at times be interested in a possible blow-up in curvature.
To this end, given $x^{m}_h \in \Vh$, we introduce the discrete curvature 
vector $\kappa^{m}_h \in \Vh$ such that
\begin{equation*} 
\int_I \pi^h [ \kappa^{m}_{h} \cdot \eta_h ]\, |x^{m}_{h,\rho}| \drho
+ \int_I \frac{x^{m}_{h,\rho} \cdot \eta_{h,\rho}}{|x^{m}_{h,\rho}|} \drho =0
\qquad \forall\ \eta_h \in \Vh.
\end{equation*}
In practice we will then monitor the quantity
\begin{equation} \label{eq:kappainv}
K^m_\infty = \max_{j=1,\ldots,J} |\kappa^m_h(q_j)|
\end{equation}
as an approximation to the maximal value of 
$|\varkappa| = \frac{|\tau_\rho|}{|x_\rho|}$. 

\subsection{Convergence experiment} \label{sec:conv3d}

\newcommand{\bgndelta}{{\hat\delta}}

We begin with a convergence experiment in order to confirm our theoretical
results from Theorem~\ref{thm:main}. To this end we
fix $\phi(p) = \sqrt{\bgndelta^2 p_1^2 + p_2^2 + p_3^2}$ 
and let $\altmu(\tau) = \frac1{\phi(\tau)}$, and then construct a suitable
right-hand side for the related flow
\[
H(y_\rho) y_t - [\Phi'(y_\rho)]_\rho = f
\]
in such a way, that the exact solution is given by
the family of self-similarly shrinking ellipses parameterized by
\begin{equation} \label{eq:y3d}
y(\rho,t)=(1-2t)^{\frac12} 
(\tfrac1{\sqrt{2}}\cos(2\pi\rho), 
\bgndelta \sin(2\pi\rho), \tfrac\bgndelta{\sqrt{2}}\cos(2\pi\rho))^T.
\end{equation}
Upon adding the correction term
\begin{equation*} 
\int_I \pi^h\left[(H(y_\rho) y_t - [\Phi'(y_\rho)]_\rho)(t_{m}) \cdot \eta_h 
\right]\! \drho 
\end{equation*}
to the right hand side of \eqref{eq:fea}, we can perform a convergence
experiment for our scheme, comparing the obtained discrete solutions with
\eqref{eq:y3d}. The results are displayed in Tables~\ref{tab:3}
and \ref{tab:4}, where we observe optimal convergence error estimates,
in line with the results proven in Theorem~\ref{thm:main}.
Here we partition the time interval $[0,T]$, with $T=0.45$, into uniform
time steps of size either $\Delta t = h$ or $\Delta t = h^2$, 
for $h = J^{-1} = 2^{-k}$, $k=6,\ldots,11$.
\begin{table}
\center
\begin{tabular}{|r|c|c|c|c|}
\hline
$J$ & $\displaystyle\max_{m=0,\ldots,M} \| y(\cdot,t_m) -  x^m_h\|_0$ & EOC
& $\displaystyle\max_{m=0,\ldots,M} \| y(\cdot,t_m) -  x^m_h\|_1$ & EOC 
\\ \hline
64  & 4.2446e-02 & ---  & 2.7051e-01 & ---  \\
128 & 2.4117e-02 & 0.82 & 1.5308e-01 & 0.82 \\
256 & 1.3110e-02 & 0.88 & 8.3054e-02 & 0.88 \\
512 & 6.8444e-03 & 0.94 & 4.3326e-02 & 0.94 \\
1024& 3.5095e-03 & 0.96 & 2.2206e-02 & 0.96 \\
2048& 1.7757e-03 & 0.98 & 1.1234e-02 & 0.98 \\
\hline
\end{tabular}
\caption{Errors for the convergence test for \eqref{eq:y3d}
over the time interval $[0,0.45]$, using $\Delta t = h$.
We also display the experimental orders of convergence (EOC).}
\label{tab:3}
\end{table}%
\begin{table}
\center
\begin{tabular}{|r|c|c|c|c|}
\hline
$J$ & $\displaystyle\max_{m=0,\ldots,M} \| y(\cdot,t_m) -  x^m_h\|_0$ & EOC
& $\displaystyle\max_{m=0,\ldots,M} \| y(\cdot,t_m) -  x^m_h\|_1$ & EOC 
\\ \hline
64  & 5.4585e-04 & ---  & 1.1777e-01 & ---  \\
128 & 1.3651e-04 & 2.00 & 5.8889e-02 & 1.00 \\
256 & 3.4129e-05 & 2.00 & 2.9445e-02 & 1.00 \\
512 & 8.5325e-06 & 2.00 & 1.4723e-02 & 1.00 \\
1024& 2.1331e-06 & 2.00 & 7.3614e-03 & 1.00 \\ 
2048& 5.3328e-07 & 2.00 & 3.6807e-03 & 1.00 \\
\hline
\end{tabular}
\caption{Errors for the convergence test for \eqref{eq:y3d}
over the time interval $[0,0.45]$, using $\Delta t = h^2$.
We also display the experimental orders of convergence (EOC).}
\label{tab:4}
\end{table}%

\subsection{Anisotropic curve shortening flow in the plane}

\begin{figure}
\center
\includegraphics[angle=-90,width=0.32\textwidth]{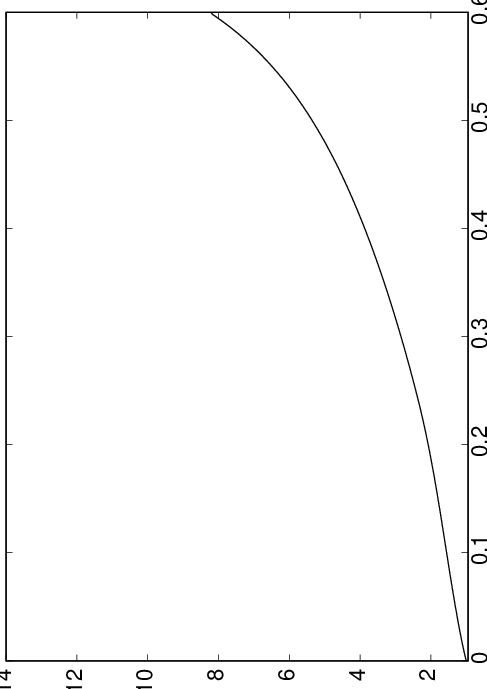}
\includegraphics[angle=-90,width=0.32\textwidth]{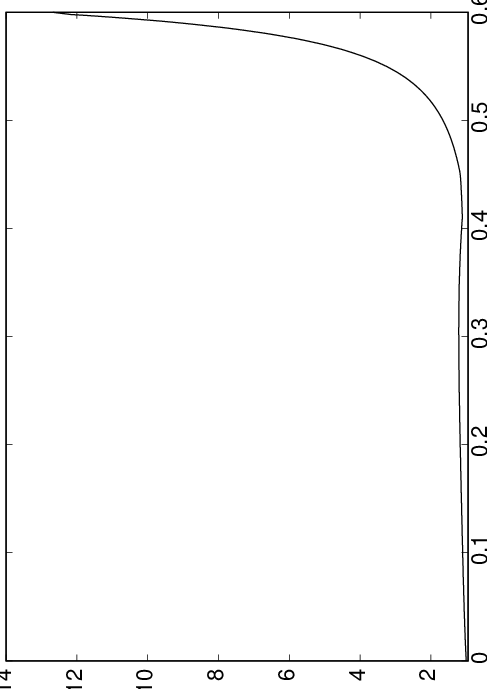}
\includegraphics[angle=-90,width=0.32\textwidth]{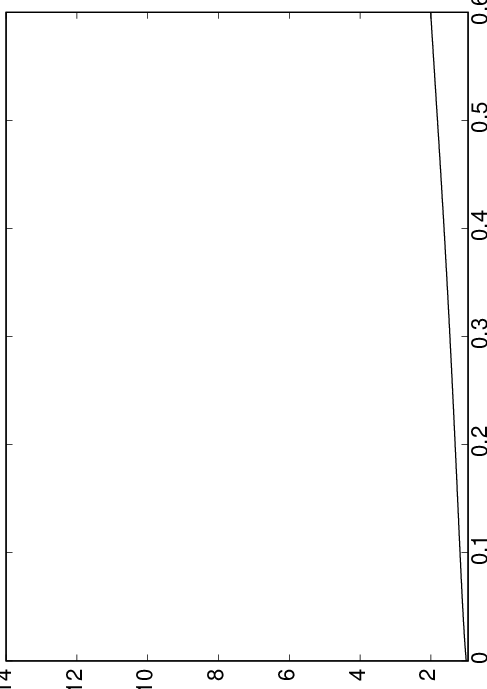}
\caption{Plots of the ratio $\ratio^m$ for \eqref{eq:fea} over time.
Left: \eqref{eq:fea}, middle: \cite{Dziuk99}, right: \cite{triplejANI}.}
\label{fig:ratios}
\end{figure}

We begin this subsection with a short investigation into the tangential
motion induced by our scheme \eqref{eq:fea}, in comparison to other schemes in
the literature. To this end, we repeat the exact same experiment from 
\cite[Figs.\ 14,15]{triplejANI}, which compared the numerical methods from
\cite{Dziuk99,triplejANI} for the anisotropic curve shortening flow 
starting from
a unit circle, for the anisotropy $\phi(p) = \sqrt{ \frac14 p_1^2 + p_2^2}$. 
In particular, we let $\altmu(\tau)=1$ and use the discretization parameters
$J=128$ and $\Delta t = 10^{-3}$. The evolution obtained from our scheme
\eqref{eq:fea} is visually indistinguishable from the results shown in
\cite[Fig.\ 14]{triplejANI}, and so we only compare plots of the ratio
$\ratio^m$ for the three schemes, see Figure~\ref{fig:ratios}. Note that the
plots shown in \cite[Fig.\ 15]{triplejANI} compared a $\phi$-weighted ratio,
for which an equidistribution result can be shown in this particular case for 
the scheme from \cite{triplejANI}. In contrast, the isotropic ratio
$\ratio^m$ shown in Figure~\ref{fig:ratios} increases for all three schemes:
the most for the scheme from \cite{Dziuk99}, the least for the scheme from
\cite{triplejANI}. For completeness we remark that repeating the experiment
with the finer discretization parameters $J=512$ and $\Delta t = 10^{-4}$ leads
to almost unchanged plots for the scheme \eqref{eq:fea} and the scheme from
\cite{triplejANI}. However, due to coalescence of mesh points, the scheme from
\cite{Dziuk99} is not able to integrate the evolution until the final time
$T=0.6$.

\begin{figure}
\center
\includegraphics[angle=-90,width=0.35\textwidth]{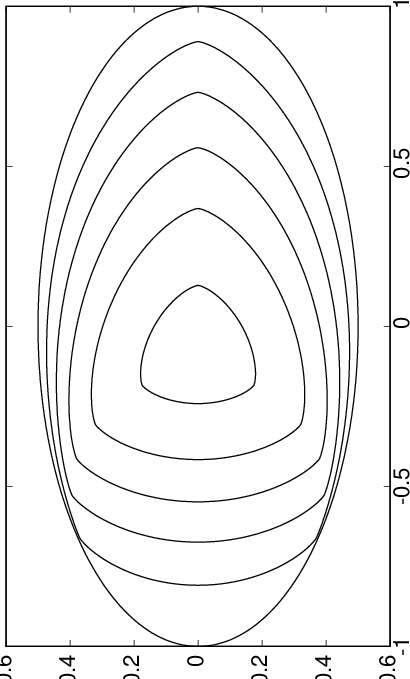}
\includegraphics[angle=-90,width=0.3\textwidth]{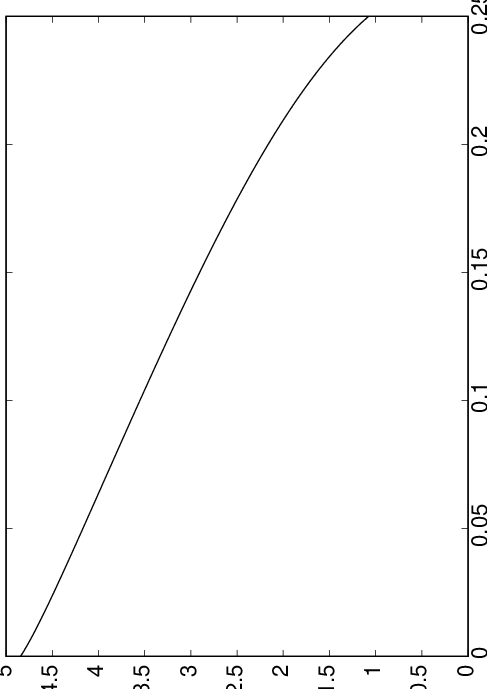}
\includegraphics[angle=-90,width=0.3\textwidth]{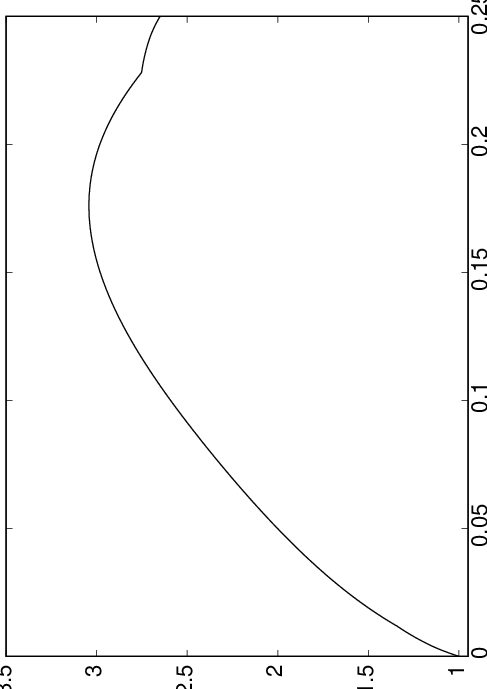}
\caption{Anisotropic curvature flow for an ellipse for the anisotropy 
\eqref{eq:gammaD99b} with $(k,\bgndelta)=(3, 0.124)$. 
Solution at times $t=0,0.05,\ldots,0.25$ on the left.
We also show plots of the discrete energy $E_\phi(x_h^m)$ (middle) and 
of the ratio $\ratio^m$ (right) over time.} 
\label{fig:ani2D99b3}
\end{figure}

For the remainder of
this subsection we again 
consider the mobility function $\altmu(\tau) = \frac1{\phi(\tau)}$. 
We begin with an anisotropy that
is not absolutely 1-homogeneous, to underline that our numerical method can
also be applied in these situations, recall Remark~\ref{rem:abshom}.
In particular, we consider an anisotropy as in \cite[(7.1)]{Dziuk99}
and \cite[(4.4a)]{fdfi}. To this end, let
\begin{equation} \label{eq:gammaD99b}
\phi( p) = | p| (1 + \bgndelta \sin(k\theta(p))),
\quad  p = | p| \binom{\cos\theta(p)}{\sin\theta(p)},\quad 
k \in \bN,\ \bgndelta \in \bR.
\end{equation}
It is not difficult to verify that this anisotropy satisfies 
\eqref{eq:sconv} if and only if $|\bgndelta| < \frac1{k^2 - 1}$.
An example computation for $k=3$ and $\bgndelta=0.124$ can be seen in
Figure~\ref{fig:ani2D99b3}, where we use as initial data an equidistributed
2:1 ellipse with unit semi major axis. During the evolution the curve
approaches as the limiting shape the so-called Wulff shape of the chosen
anisotropy. Here we recall from \cite{FonsecaM91} that the Wulff shape
is the solution of the isoperimetric problem for 
the energy \eqref{eq:Ephi} in the plane.

\begin{figure}
\center
\mbox{
\includegraphics[angle=-90,width=0.24\textwidth]{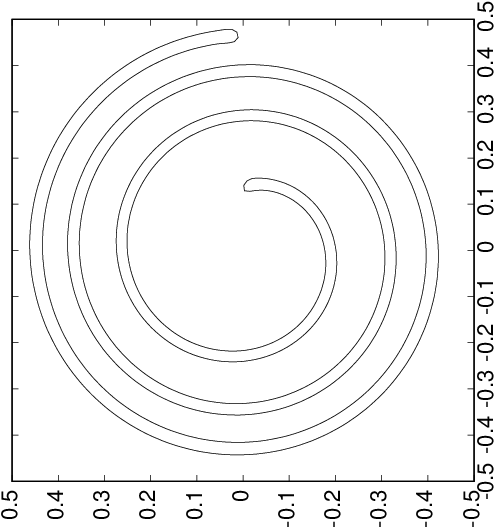}
\includegraphics[angle=-90,width=0.24\textwidth]{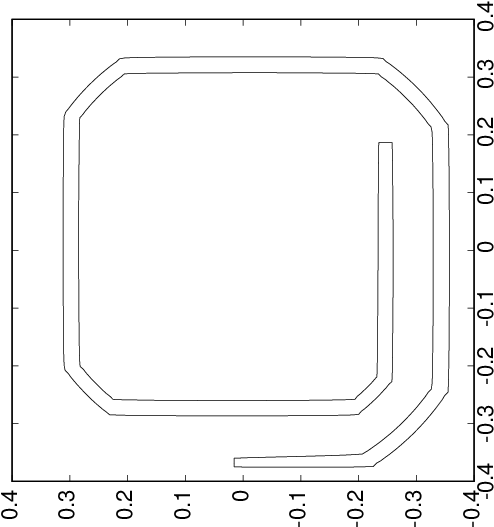}
\includegraphics[angle=-90,width=0.24\textwidth]{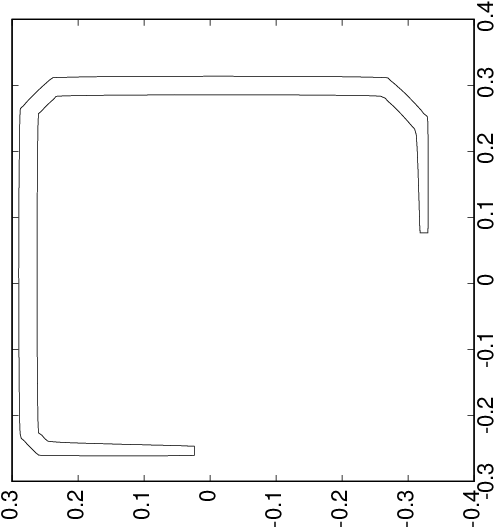}
\includegraphics[angle=-90,width=0.24\textwidth]{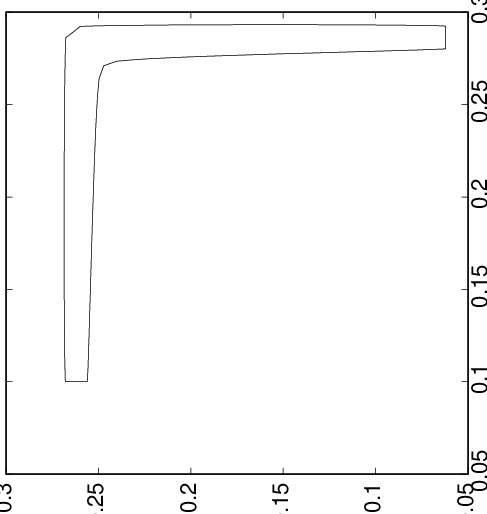}
}
\caption{Anisotropic curvature flow for a spiral for the anisotropy 
\eqref{eq:bgnL2}. 
Solution at times $t=0,0.01,0.015,0.019$.
}
\label{fig:ani2bgnL2d001}
\end{figure}

With our next numerical experiment we would like to demonstrate that our
scheme can also deal with nearly crystalline anisotropies. To this end, we
choose as anisotropy the regularized $\ell^1$--norm
\begin{equation} \label{eq:bgnL2}
 \phi(p) = \sum_{i = 1}^d\sqrt{(1-\bgndelta^2) p_i^2 + \bgndelta^2|p|^2}, \quad 
\bgndelta = 0.01. 
\end{equation}
A simulation starting from a spiral shape is shown in
Figure~\ref{fig:ani2bgnL2d001}, 
where we notice the developing facets and corners
due to the crystalline nature of the chosen anisotropy.
For this experiment we use the finer discretization parameters
$J=1024$ and $\Delta t = 5 \times 10^{-8}$.

\subsection{Isotropic curve shortening flow in $\bR^3$}

\begin{figure}
\center
\includegraphics[angle=-0,width=0.43\textwidth]{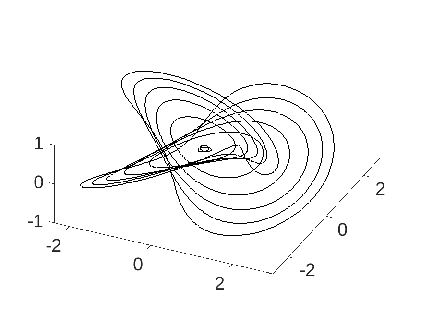}
\includegraphics[angle=-0,width=0.43\textwidth]{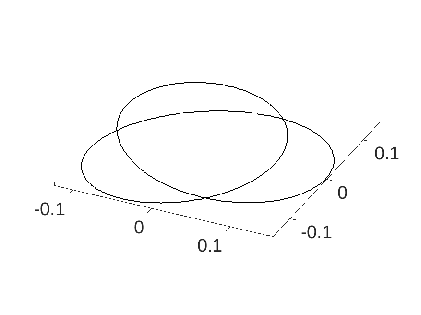}
\includegraphics[angle=-90,width=0.3\textwidth]{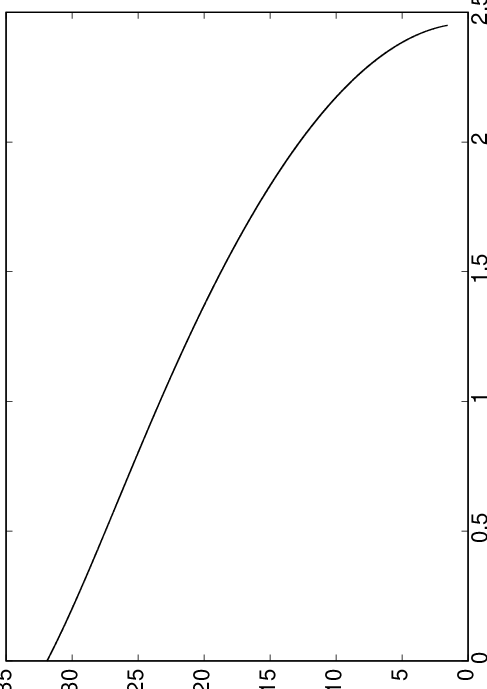}
\includegraphics[angle=-90,width=0.3\textwidth]{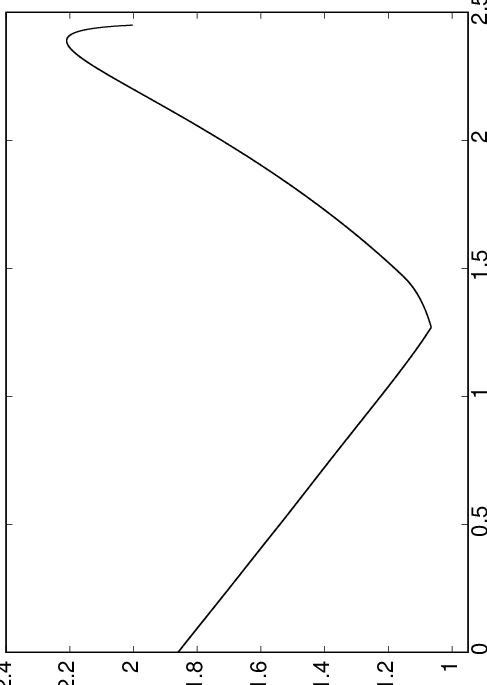}
\includegraphics[angle=-90,width=0.3\textwidth]{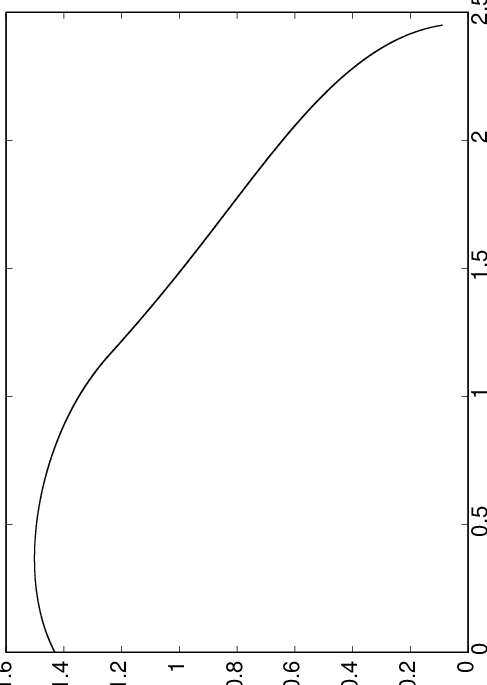}
\caption{Isotropic curve shortening flow for the trefoil knot
\eqref{eq:trefoil}. On the left,
$x_h^m$ at times $t=0,0.5,\ldots,2,T=2.45$, with $x_h^M$ on the right.
Below we show plots of $E_\phi(x_h^m)$, $\ratio^m$ and $1/K^m_\infty$ over 
time.}
\label{fig:knot}
\end{figure}

The remainder of our numerical experiments are for the two-codimensional case,
and from now on we always choose the constant mobility $\altmu(\tau) = 1$.
In this particular subsection we in addition consider the isotropic case,
$\phi(p)=|p|$.
The first experiment is for a trefoil knot in $\bR^3$, and in particular
the initial curve is given by
\begin{equation}
x_0(\rho) = ((2+\cos(6\pi\rho))\cos(4\pi\rho),
  (2+\cos(6\pi\rho))\sin(4\pi\rho), \sin(6\pi\rho) )^T,
 \quad \rho \in I.
\label{eq:trefoil}
\end{equation}
See Figure~\ref{fig:knot} for the numerical results, which agree very well with
the results from \cite[Figure~1]{curves3d}. Observe that the knot
approaches a double covering of a circle within a hyperplane of $\bR^3$.

\begin{figure}
\center
\includegraphics[angle=-0,width=0.3\textwidth]{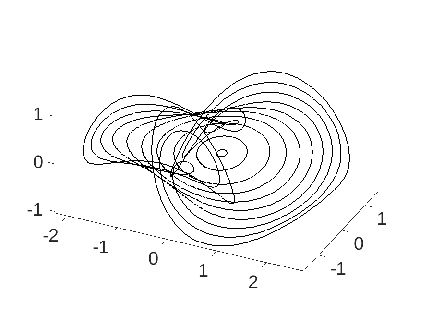}
\includegraphics[angle=-0,width=0.3\textwidth]{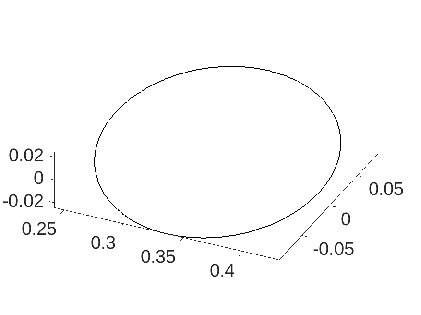}
\includegraphics[angle=-0,width=0.3\textwidth]{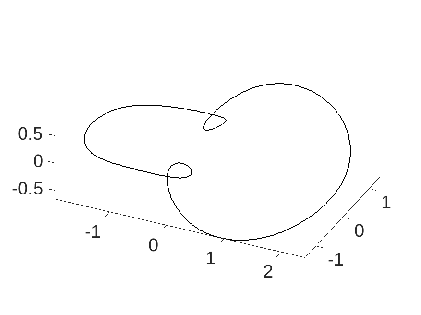}
\includegraphics[angle=-90,width=0.3\textwidth]{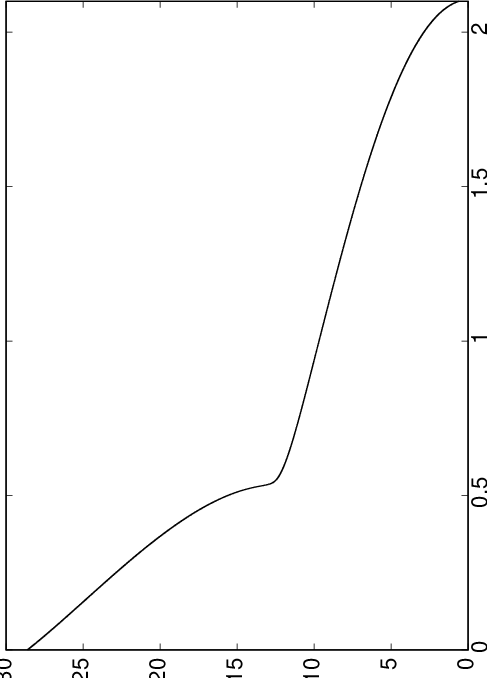}
\includegraphics[angle=-90,width=0.3\textwidth]{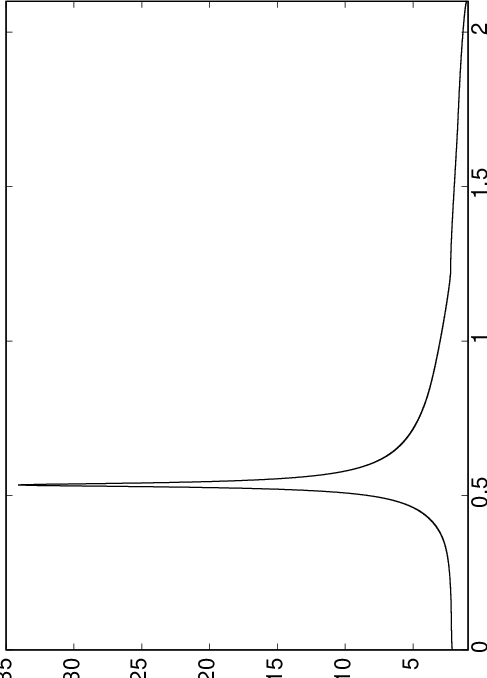}
\includegraphics[angle=-90,width=0.3\textwidth]{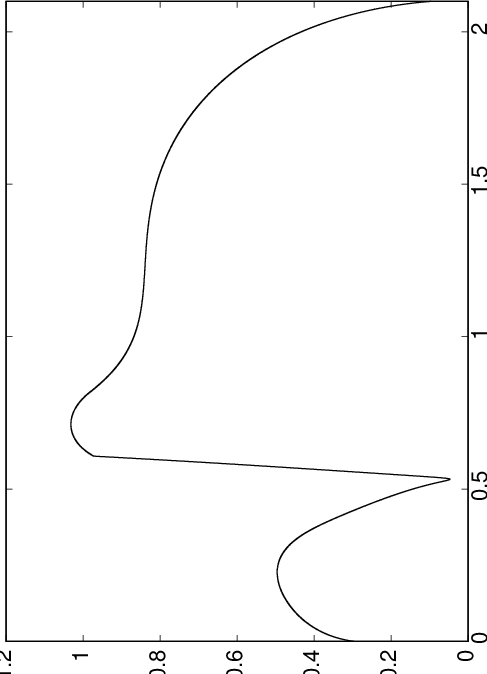}
\caption{Isotropic curve shortening flow for the two interlocked rings
\eqref{eq:irings}. On the left,
$x_h^m$ at times $t=0,0.25,\ldots,2,T=2.1$, in the middle $x^M_h$,
with $x_h^m$ at time $t=0.5$ on the right.
Below we also show plots of $E_\phi(x_h^m)$, 
$\ratio^m$ and $1/K^m_\infty$ over time.}
\label{fig:irings}
\end{figure}

The next experiment is for two interlocked rings in $\bR^3$, and in particular
the initial curve is given by
\begin{equation}
x_0(\rho) = \tfrac18 
\begin{pmatrix}
10 (\cos(2\pi\rho)+\cos(6\pi\rho))+\cos(4\pi\rho)
+\cos(8\pi\rho) \\
6\sin(2\pi\rho)+10\sin(6\pi\rho) \\
4\sin(6\pi\rho)\sin(5\pi\rho)+4\sin(8\pi\rho)-2\sin(12\pi\rho)
\end{pmatrix}
\quad \rho \in I.
\label{eq:irings}
\end{equation}
See Figure~\ref{fig:irings} for the numerical results, where we can observe a
singularity in the flow. The curve forms two loops that shrink, developing
into two dove-tails. Here the continuous problem develops a singularity, with
the curvature blowing up. From the plot of the inverse
of the magnitude of the discrete curvature, \eqref{eq:kappainv}, 
we can clearly see the effect of the singularity also on the discrete level. 
Let us emphasize that existence and error estimates for the discrete solution
only hold as long as the required smoothness assumptions are satisfied, i.e.\
before the formation of the singularity.
Our numerical scheme, however, simply integrates
through the singularity and eventually approaches a circle that shrinks to a
point within a hyperplane. Notions of weak solutions that allow an extension of the solution
beyond the singularity have been proposed in \cite{Angenent91,Deckelnick97} for the isotropic case, but
it is very difficult to obtain a convergence result after the onset of singularities.

\begin{figure}
\center
\includegraphics[angle=-0,width=0.3\textwidth]{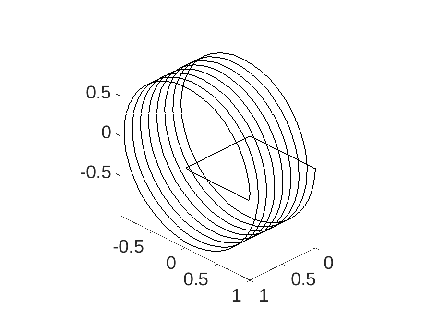}
\includegraphics[angle=-0,width=0.3\textwidth]{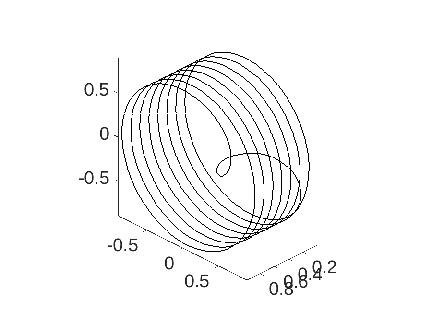}
\includegraphics[angle=-0,width=0.3\textwidth]{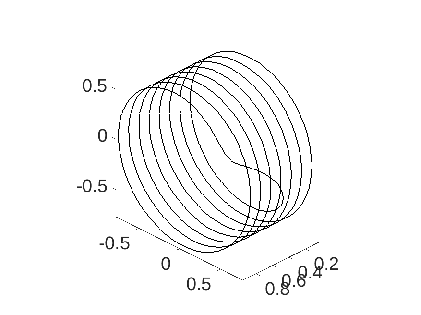}
\includegraphics[angle=-0,width=0.3\textwidth]{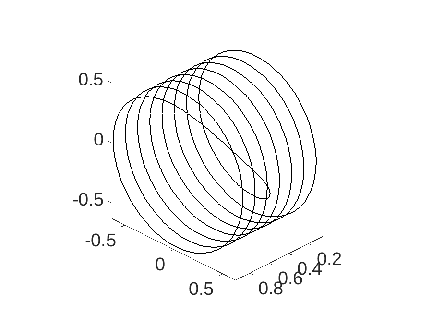}
\includegraphics[angle=-0,width=0.3\textwidth]{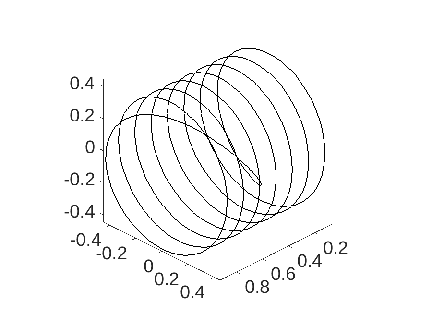}
\includegraphics[angle=-0,width=0.3\textwidth]{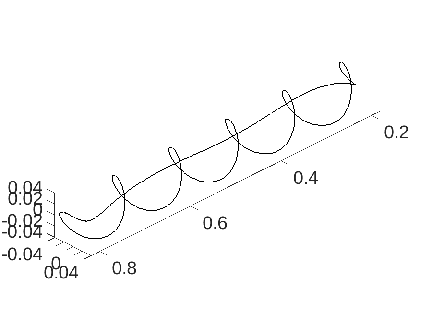}
\caption{Isotropic curve shortening flow for the helix \eqref{eq:helix}. 
We show $x_h^m$ at times $t=0,0.1,\ldots,T=0.5$.
}
\label{fig:helix}
\end{figure}

In our final experiment for the isotropic setting we consider
a closed helix in $\bR^3$, as in \cite[Figure~2]{curves3d}. 
Here the open helix is defined by
\begin{equation}
 x_0(\varrho) = (\sin(16\,\pi\varrho), \cos(16\,\pi\,\varrho), \varrho)^T
\,,\quad \varrho \in [0,1]\,,
\label{eq:helix}
\end{equation}
and the initial curve is constructed from \eqref{eq:helix}
by connecting $x_0(0)$ and $x_0(1)$ with a polygon that visits the origin
and $(0,0,1)^T$. The evolution of the helix under curve shortening flow can
be seen in Figure~\ref{fig:helix}.

\subsection{Anisotropic curve shortening flow in $\bR^3$}

\begin{figure}
\center
\includegraphics[angle=-0,width=0.43\textwidth]{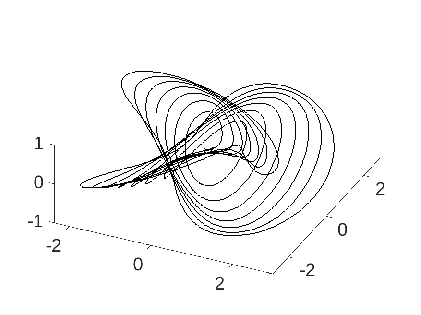}
\includegraphics[angle=-0,width=0.43\textwidth]{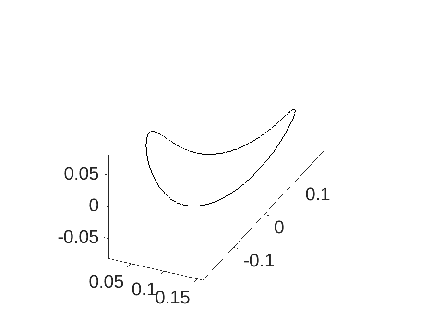}
\caption{Anisotropic curve shortening flow for the trefoil knot
\eqref{eq:trefoil} for the anisotropy 
$\phi(p) = \sqrt{p_1^2 + \tfrac14 (p_2^2 + p_3^2)}$. On the left,
$x_h^m$ at times $t=0,0.5,\ldots,T=3.5$, with $x_h^M$ on the right.
}
\label{fig:aniknot}
\end{figure}

\begin{figure}
\center
\includegraphics[angle=-0,width=0.3\textwidth]{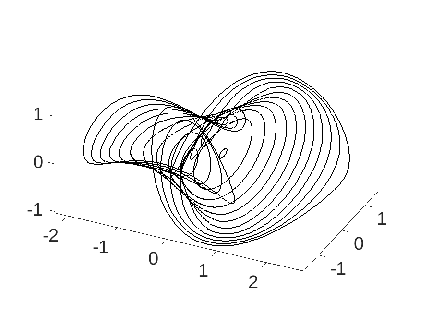}
\includegraphics[angle=-0,width=0.3\textwidth]{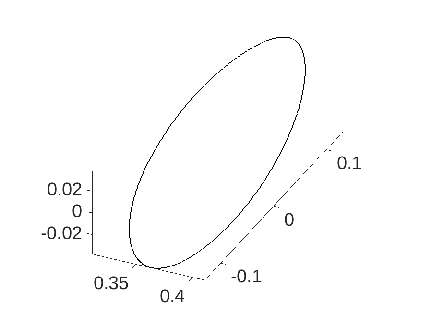}
\includegraphics[angle=-0,width=0.3\textwidth]{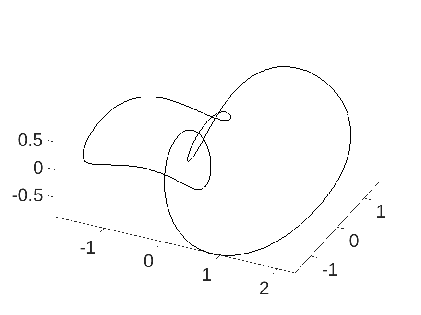}
\includegraphics[angle=-90,width=0.3\textwidth]{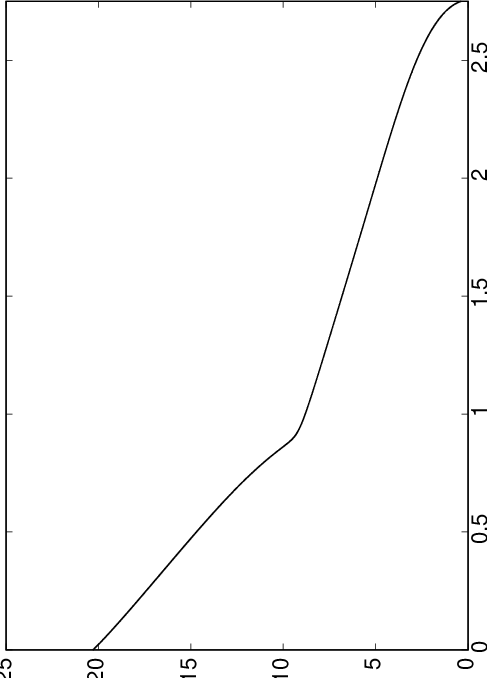}
\includegraphics[angle=-90,width=0.3\textwidth]{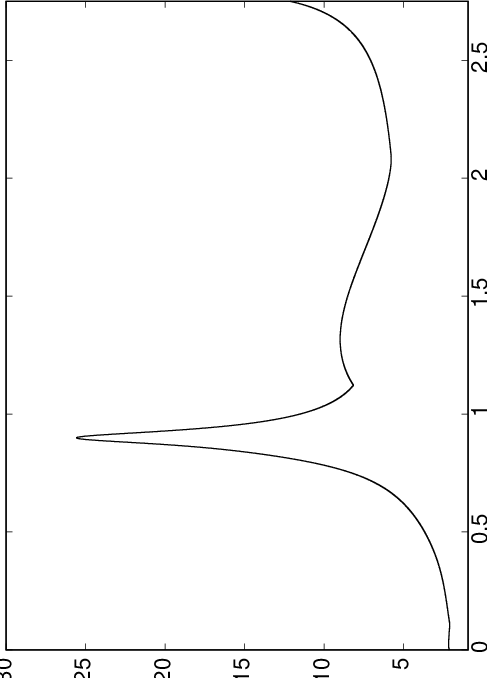}
\includegraphics[angle=-90,width=0.3\textwidth]{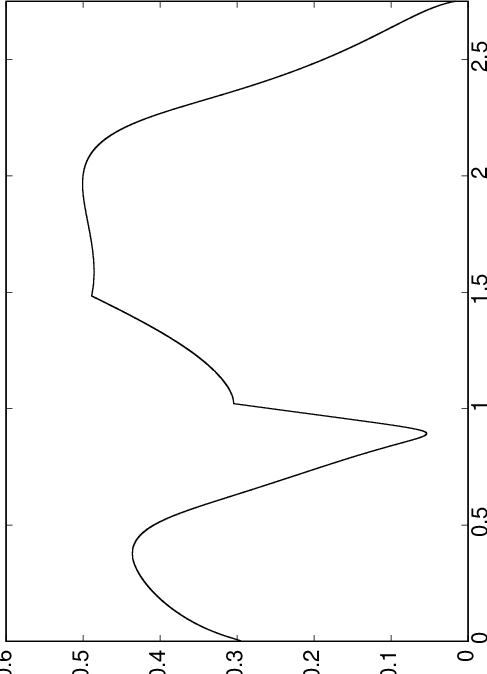}
\caption{Anisotropic curve shortening flow for the two interlocked rings
\eqref{eq:irings} for the anisotropy 
$\phi(p) = \sqrt{p_1^2 + \tfrac14 (p_2^2 + p_3^2)}$. On the left,
$x_h^m$ at times $t=0,0.25,\ldots,2.5,T=2.75$, in the middle $x_h^M$,
with $x_h^m$ at time $t=0.5$ on the right.
Below we also show plots of $E_\phi(x_h^m)$, 
$\ratio^m$ and $1/K^m_\infty$ over time.}
\label{fig:aniirings}
\end{figure}

\begin{figure}
\center
\includegraphics[angle=-0,width=0.3\textwidth]{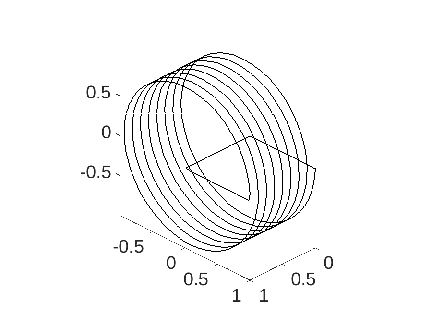}
\includegraphics[angle=-0,width=0.3\textwidth]{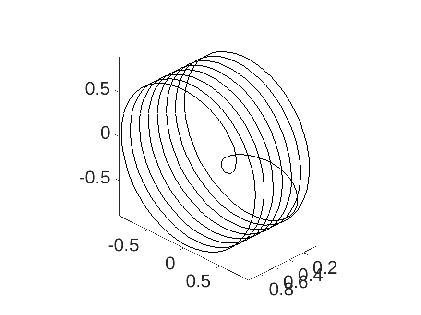}
\includegraphics[angle=-0,width=0.3\textwidth]{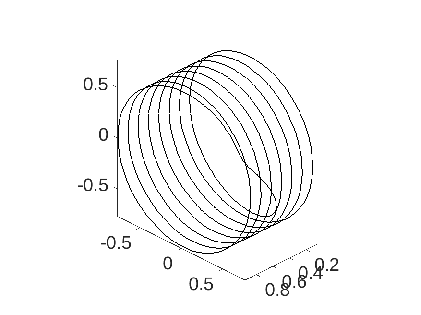}
\includegraphics[angle=-0,width=0.3\textwidth]{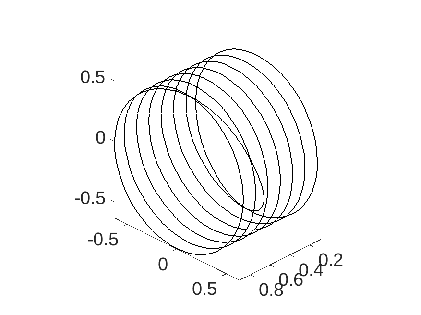}
\includegraphics[angle=-0,width=0.3\textwidth]{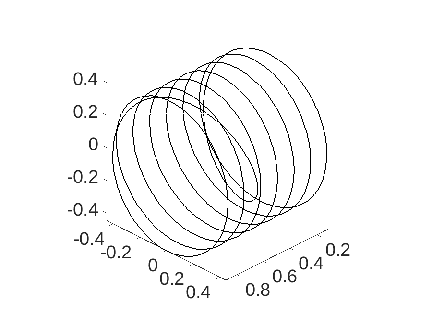}
\includegraphics[angle=-0,width=0.3\textwidth]{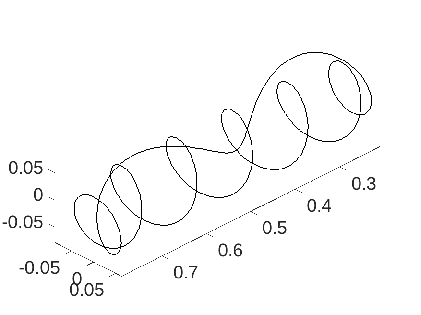}
\caption{Anisotropic curve shortening flow for the helix \eqref{eq:helix} for 
the anisotropy $\phi(p) = \sqrt{p_1^2 + \tfrac14 (p_2^2 + p_3^2)}$. 
We show $x_h^m$ at times $t=0,0.2,\ldots,T=1$.
}
\label{fig:anihelix}
\end{figure}

\begin{figure}
\center
\includegraphics[angle=-0,width=0.3\textwidth]{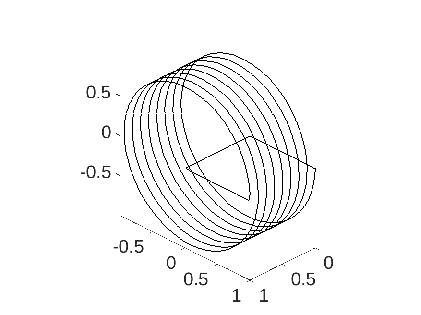}
\includegraphics[angle=-0,width=0.3\textwidth]{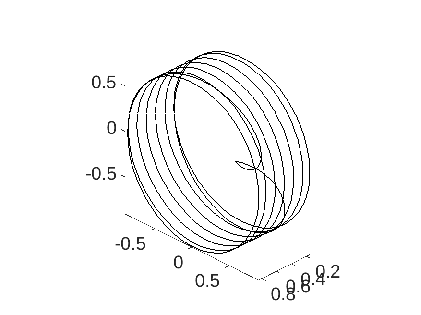}
\includegraphics[angle=-0,width=0.3\textwidth]{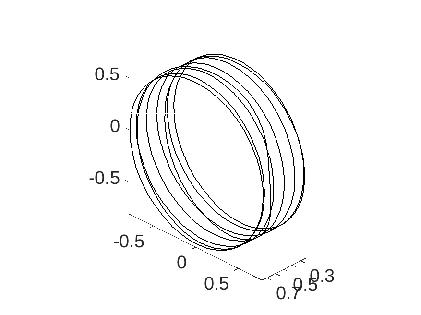}
\includegraphics[angle=-0,width=0.3\textwidth]{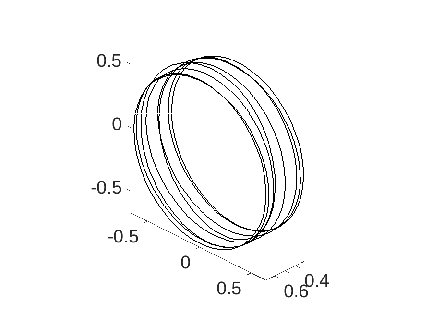}
\includegraphics[angle=-0,width=0.3\textwidth]{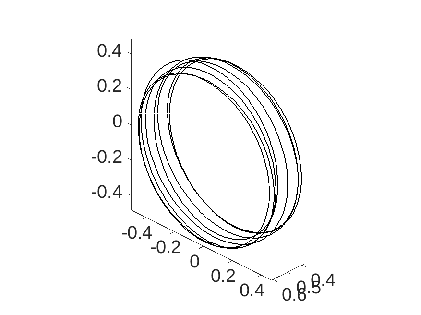}
\includegraphics[angle=-0,width=0.3\textwidth]{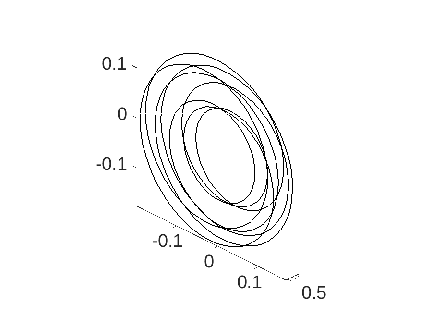}
\includegraphics[angle=-90,width=0.3\textwidth]{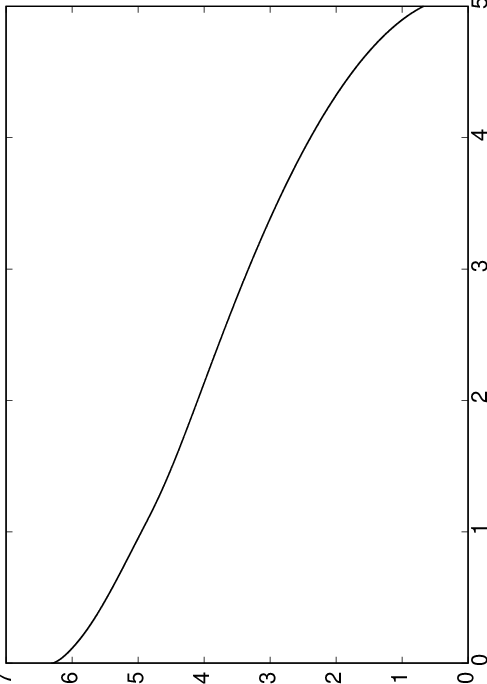}
\includegraphics[angle=-90,width=0.3\textwidth]{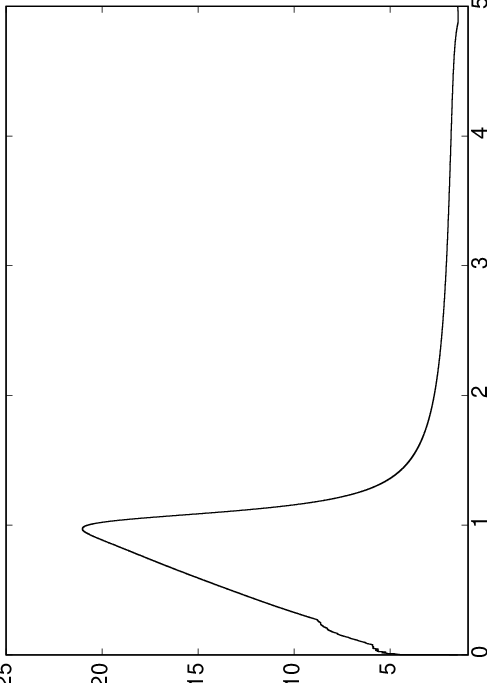}
\includegraphics[angle=-90,width=0.3\textwidth]{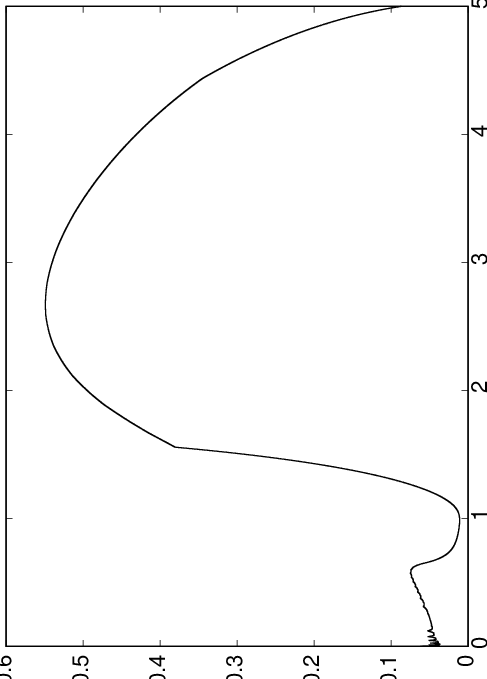}
\caption{Anisotropic curve shortening flow for the helix \eqref{eq:helix} for 
the anisotropy $\phi(p) = \sqrt{p_1^2 + \tfrac14 (p_2^2 + p_3^2)}$. 
We show $x_h^m$ at times $t=0,1,\ldots,T=5$.
Below we also show plots of $E_\phi(x_h^m)$, 
$\ratio^m$ and $1/K^m_\infty$ over time.}
\label{fig:ani01helix}
\end{figure}

Still for the constant mobility $\altmu(\tau) = 1$, we now repeat the
experiments from the previous subsection, but now for the anisotropy
\[
\phi(p) = \sqrt{p_1^2 + \bgndelta^2 (p_2^2 + p_3^2)}, \quad \bgndelta = 0.5.
\]
This means that for the
evolving curves it will be energetically favourable to have tangents in the
$y-z$ plane, which should result in the curve itself trying to migrate into
that hyperplane.
See the results in Figures~\ref{fig:aniknot}, \ref{fig:aniirings}
and \ref{fig:anihelix}.
We observe that in Figure~\ref{fig:aniirings} it is no longer clear that a
singularity occurs for the continuous problem. In particular, the evolution
of the inverse of the maximal discrete curvature appears to indicate
that the curvature for the continuous problem remains bounded.
For the evolution of the initial helix there seems to be little qualitative
difference compared to the isotropic results.
However, on repeating the simulation for the stronger anisotropy 
$\phi(p) = \sqrt{p_1^2 + 0.01 (p_2^2 + p_3^2)}$ yields the evolution in
Figure~\ref{fig:ani01helix}, which clearly shows the effect of the anisotropy.
In particular, it can be seen that the curve attempts to avoid any tangent
vectors that have a non-zero $x$-component, as is to be expected from
the chosen anisotropy. We also note the large values of
$K^m_\infty$ around time $t=1$, which may indicate the development of a
possible singularity in the flow.

\begin{figure}
\center
\includegraphics[angle=-0,width=0.3\textwidth]{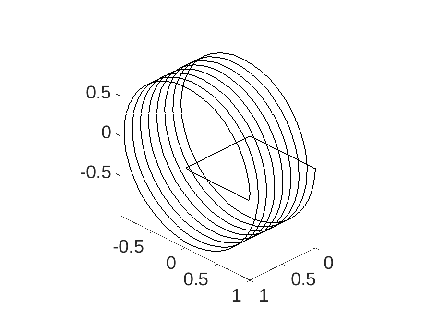}
\includegraphics[angle=-0,width=0.3\textwidth]{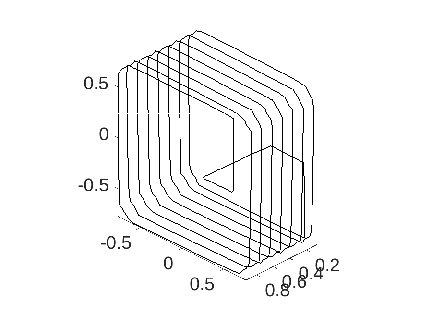}
\includegraphics[angle=-0,width=0.3\textwidth]{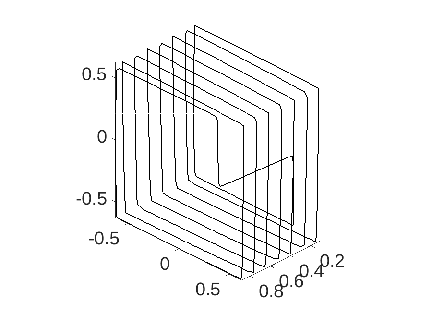}
\includegraphics[angle=-0,width=0.3\textwidth]{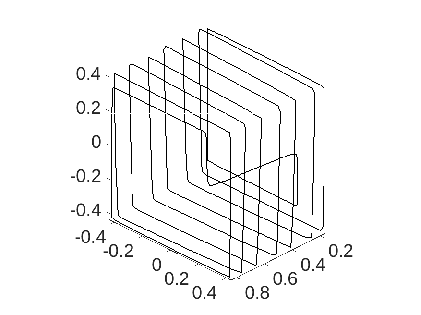}
\includegraphics[angle=-0,width=0.3\textwidth]{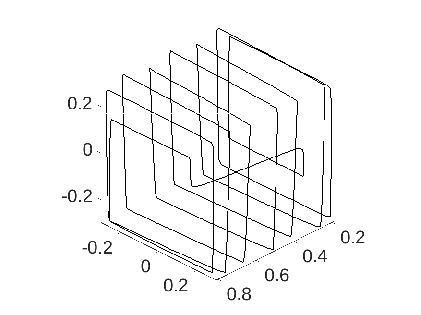}
\includegraphics[angle=-0,width=0.3\textwidth]{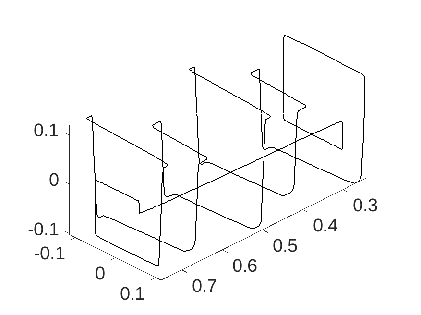}
\caption{Anisotropic curve shortening flow for the helix \eqref{eq:helix} for 
the anisotropy \eqref{eq:bgnL2}.
We show $x_h^m$ at times $t=0,0.1,0.2,0.3,0.25,T=0.39$.
}
\label{fig:aniL3001helix}
\end{figure}

We end this section with another numerical simulation 
for the regularized $\ell^1$--norm \eqref{eq:bgnL2}.
This nearly crystalline anisotropy forces the curve to have tangents aligned
with the three coordinate axes. Starting from an initial helix as before,
we can observe that effect in
Figure~\ref{fig:aniL3001helix}. For this experiment we used the smaller time
step size $\Delta t = 10^{-6}$.

\clearpage
\begin{appendix}
\renewcommand{\theequation}{\Alph{section}.\arabic{equation}}
\setcounter{equation}{0}
\section{Appendix} \label{sec:appA}

{\it Proof of Lemma~\ref{lem:nonlinearinterpol}}: 
Let $y$ be given as in Lemma~\ref{lem:nonlinearinterpol}.
We define $A \in [L^\infty(I)]^{d \times d}$ by 
$A_{| I_j} = \frac{1}{| I_j|} \int_{I_j} \Phi''(y_\rho) \drho$, 
$j=1,\ldots J$. Clearly it holds that
\begin{equation} \label{eq:difaj}
\| A - \Phi''(y_\rho) \|_{0,\infty} \leq h | \Phi''(y_\rho) |_{1,\infty} \leq C h \Vert y \Vert_{2,\infty} \leq CC_1 h. 
\end{equation}
Moreover, similarly to \eqref{eq:convex3} there exists $\sigma_1 > 0$ 
independent of $y$ such that
\begin{equation} \label{eq:Alb}
A \xi \cdot \xi \geq \sigma_1 |\xi|^2 \quad \forall\ \xi \in \bR^d \quad
\text{in}\ I.
\end{equation}
Observing that
\begin{align*}
\int_I A y_\rho \cdot \eta_{h,\rho} \drho &
= \sum_{j=1}^J A_{| I_j} \eta_{h,\rho} \cdot \int_{I_j} y_\rho \drho 
= \sum_{j=1}^J A_{| I_j} \eta_{h,\rho} \cdot \int_{I_j} (\pi^h y)_\rho \drho 
\\ &
= \int_I A (\pi^h y)_\rho \cdot \eta_{h,\rho} \drho \quad
\forall\ \eta_h \in \Vh,
\end{align*}
we may write \eqref{eq:defhatx} in the form
\begin{align*}
& \int_I A (Q_h y -  \pi^h y)_\rho \cdot \eta_{h,\rho} \drho 
+ \int_I (Q_h y - \pi^h y) \cdot \eta_h \drho \\ & \qquad
= \int_I (A - \Phi''(y_\rho)) (Q_h y - y)_\rho \cdot \eta_{h,\rho} \drho 
- \int_I (\pi^h y - y) \cdot \eta_h \drho 
\quad \forall\ \eta_h \in \Vh.
\end{align*}
On testing the above relation with $\eta_h = Q_h y - \pi^h y$, 
and recalling \eqref{eq:Alb}, \eqref{eq:difaj} and \eqref{eq:estpih}, 
we infer that
\begin{align*}
& \min\{\sigma_1,1\} \| Q_h y - \pi^h y \|_1^2
\leq C h | Q_h y - y |_1 | Q_h y - \pi^h y |_1 
+ \| \pi^h y - y \|_0 \| Q_h y - \pi^h y \|_0 \\ & \qquad
\leq C h (| Q_h y - \pi^h y |_1 + | \pi^h y - y|_1 ) | Q_h y - \pi^h y |_1 
+ C h^2 |y|_2 \| Q_h y - \pi^h y \|_0 \\ & \qquad
\leq C h | Q_h y - \pi^h y |_1^2 + C h^2 |y|_2 \| Q_h y - \pi^h y \|_1
\leq \tfrac12 \min\{\sigma_1,1\} \| Q_h y - \pi^h y \|_1^2 + C h^4,
\end{align*}
provided that $0<h \leq h_1$ and $h_1>0$ is sufficiently small. Hence we have that
\begin{equation} \label{eq:appAPhpih}
\| Q_h y - \pi^h y \|_1 \leq C h^2.
\end{equation}
The estimate \eqref{eq:approx} now follows with the help of 
\eqref{eq:estpih}, while \eqref{eq:estpih}, \eqref{eq:inverse} and \eqref{eq:appAPhpih} imply that
\[
| y - Q_h y |_{1,\infty} \leq | y - \pi^h y |_{1,\infty} + | \pi^h y - Q_h y |_{1,\infty} \leq C h | y |_{2, \infty} + Ch^{-\frac12} | \pi^h y - Q_h y |_1 \leq Ch.
\]
In particular, $|(Q_h y)_\rho| \geq | y_\rho | - | y-Q_h y|_{1,\infty} \geq c_1 -Ch \geq \tfrac12 c_1$ provided that $h_1$ is small enough, and similarly
$| (Q_h y)_\rho | \leq 2 C_1$. Hence we have shown \eqref{eq:approxinf}.
 
Let us now consider the case where $y$ depends in addition on time.
Differentiating \eqref{eq:defhatx} with respect to $t$ yields that 
\begin{align} 
& \int_I \Phi''(y_\rho) (Q_h y - y)_{t,\rho} \cdot \eta_{h,\rho} \drho
+ \int_I (Q_h y - y)_t \cdot \eta_h \drho \nonumber \\ & \qquad
+ \int_I \Phi'''(y_\rho) (y_{t,\rho}, (Q_h y-y)_\rho, \eta_{h,\rho}) \drho= 0
\quad \forall\ \eta_h \in \Vh.  \label{eq:defhatxdt}
\end{align}
Let $A \in L^\infty(0,T;[L^\infty(I)]^{d \times d})$ be defined analogously
to above, and $B \in L^\infty(0,T;[L^\infty(I)]^{d \times d})$ via
$B_{| I_j} =\frac{1}{|I_j|} \int_{I_j} \Phi'''(y_\rho)
(y_{t,\rho},\cdot,\cdot)\drho$ 
on $[0,T]$. Then we can write \eqref{eq:defhatxdt} in the form
\begin{align*}
& \int_I A (Q_h y -  \pi^h y)_{t,\rho} \cdot \eta_{h,\rho} \drho 
+ \int_I (Q_h y - \pi^h y)_t \cdot \eta_h \drho \nonumber \\ & \quad
=  \int_I (A - \Phi''(y_\rho)) (Q_h y - y)_{t,\rho} \cdot \eta_{h,\rho} \drho 
- \int_I (\pi^h y - y)_t \cdot \eta_h \drho \\ & \quad 
+ \int_I B (Q_h y - \pi^h y)_\rho \cdot \eta_{h,\rho} \drho
+ \int_I (B - \Phi'''(y_\rho) (y_{t,\rho},\cdot,\cdot)) 
(Q_h y -  y)_\rho \cdot \eta_{h,\rho} \drho  \quad \forall\ \eta_h \in \Vh,
\end{align*}
where we have used the relation $\int_I B y_\rho \cdot \eta_{h,\rho} \drho = \int_I B (\pi^h y)_\rho \cdot \eta_{h,\rho} \drho$. If
we choose $\eta_h= (Q_h y - \pi^h y)_t \in \Vh$, use \eqref{eq:appAPhpih}, 
the bound $\| B - \Phi'''(y_\rho) (y_{t,\rho},\cdot,\cdot) \|_{0,\infty} \leq CC_1 h$ and argue
similarly as above we obtain
\begin{displaymath}
\| (Q_h y - \pi^h y)_t \|_1 \leq Ch^2,
\end{displaymath}
from which we infer the first estimate in \eqref{eq:approxtime}. Combining this bound with \eqref{eq:inverse} and \eqref{eq:estpih} we finally obtain
\begin{align*}
| (Q_h y)_t |_{1,\infty} &  \leq  | (Q_h y - \pi^h y)_t |_{1,\infty} + | (\pi^h y -y)_t |_{1,\infty} + | y_t | _{1,\infty}  \\
& \leq C h^{-\frac{1}{2}} | (Q_h y - \pi^h y)_t |_1 + Ch |  y_t |_{2,\infty} + C_1 \leq C,
\end{align*}
which completes the proof of Lemma \ref{lem:nonlinearinterpol}.

\end{appendix}

\def\soft#1{\leavevmode\setbox0=\hbox{h}\dimen7=\ht0\advance \dimen7
  by-1ex\relax\if t#1\relax\rlap{\raise.6\dimen7
  \hbox{\kern.3ex\char'47}}#1\relax\else\if T#1\relax
  \rlap{\raise.5\dimen7\hbox{\kern1.3ex\char'47}}#1\relax \else\if
  d#1\relax\rlap{\raise.5\dimen7\hbox{\kern.9ex \char'47}}#1\relax\else\if
  D#1\relax\rlap{\raise.5\dimen7 \hbox{\kern1.4ex\char'47}}#1\relax\else\if
  l#1\relax \rlap{\raise.5\dimen7\hbox{\kern.4ex\char'47}}#1\relax \else\if
  L#1\relax\rlap{\raise.5\dimen7\hbox{\kern.7ex
  \char'47}}#1\relax\else\message{accent \string\soft \space #1 not
  defined!}#1\relax\fi\fi\fi\fi\fi\fi}


\begin{thebibliography}{10}

\bibitem{AmbrosioS96}
{\sc L.~Ambrosio and H.~M. Soner}, {\em Level set approach to mean curvature
  flow in arbitrary codimension}, J. Differential Geom., 43 (1996),
  pp.~693--737.

\bibitem{Angenent91}
{\sc S.~Angenent}, {\em Parabolic equations for curves on surfaces. {II}.
  {I}ntersections, blow-up and generalized solutions}, Ann. of Math., 133
  (1991), pp.~171--215.

\bibitem{BanschDGP23}
{\sc E.~B\"{a}nsch, K.~Deckelnick, H.~Garcke, and P.~Pozzi}, {\em Interfaces:
  Modeling, Analysis, Numerics}, vol.~51 of Oberwolfach Seminars,
  Birkh\"{a}user/Springer, Cham, 2023.

\bibitem{triplejANI}
{\sc J.~W. Barrett, H.~Garcke, and R.~N\"urnberg}, {\em Numerical approximation
  of anisotropic geometric evolution equations in the plane}, IMA J. Numer.
  Anal., 28 (2008), pp.~292--330.

\bibitem{curves3d}
\leavevmode\vrule height 2pt depth -1.6pt width 23pt, {\em Numerical
  approximation of gradient flows for closed curves in {${\mathbb R}^d$}}, IMA
  J. Numer. Anal., 30 (2010), pp.~4--60.

\bibitem{fdfi}
\leavevmode\vrule height 2pt depth -1.6pt width 23pt, {\em The approximation of
  planar curve evolutions by stable fully implicit finite element schemes that
  equidistribute}, Numer. Methods Partial Differential Equations, 27 (2011),
  pp.~1--30.

\bibitem{bgnreview}
\leavevmode\vrule height 2pt depth -1.6pt width 23pt, {\em Parametric finite
  element approximations of curvature driven interface evolutions}, in Handb.
  Numer. Anal., A.~Bonito and R.~H. Nochetto, eds., vol.~21, Elsevier,
  Amsterdam, 2020, pp.~275--423.

\bibitem{BellettiniP96}
{\sc G.~Bellettini and M.~Paolini}, {\em Anisotropic motion by mean curvature
  in the context of {F}insler geometry}, Hokkaido Math. J., 25 (1996),
  pp.~537--566.

\bibitem{BenesKS22}
{\sc M.~Bene\v{s}, M.~Kol\'{a}\v{r}, and D.~\v{S}ev\v{c}ovi\v{c}}, {\em
  Qualitative and numerical aspects of a motion of a family of interacting
  curves in space}, SIAM J. Appl. Math., 82 (2022), pp.~549--575.

\bibitem{BinzK21preprint}
{\sc T.~Binz and B.~Kov{\'a}cs}, {\em A convergent finite element algorithm for
  mean curvature flow in arbitrary codimension}.
\newblock arXiv:2107.10577, 2021.

\bibitem{CarliniFF07}
{\sc E.~Carlini, M.~Falcone, and R.~Ferretti}, {\em A semi-{L}agrangian scheme
  for the curve shortening flow in codimension-2}, J. Comput. Phys., 225
  (2007), pp.~1388--1408.

\bibitem{Davis04}
{\sc T.~A. Davis}, {\em Algorithm 832: {UMFPACK} {V}4.3---an
  unsymmetric-pattern multifrontal method}, ACM Trans. Math. Software, 30
  (2004), pp.~196--199.

\bibitem{Deckelnick97}
{\sc K.~Deckelnick}, {\em Weak solutions of the curve shortening flow}, Calc.
  Var. Partial Differential Equations, 5 (1997), pp.~489--510.

\bibitem{DeckelnickD95}
{\sc K.~Deckelnick and G.~Dziuk}, {\em On the approximation of the curve
  shortening flow}, in Calculus of Variations, Applications and Computations
  (Pont-\`a-Mousson, 1994), C.~Bandle, J.~Bemelmans, M.~Chipot, J.~S.~J.
  Paulin, and I.~Shafrir, eds., vol.~326 of Pitman Res. Notes Math. Ser.,
  Longman Sci. Tech., Harlow, 1995, pp.~100--108.

\bibitem{DeckelnickDE05}
{\sc K.~Deckelnick, G.~Dziuk, and C.~M. Elliott}, {\em Computation of geometric
  partial differential equations and mean curvature flow}, Acta Numer., 14
  (2005), pp.~139--232.

\bibitem{finsler}
{\sc K.~Deckelnick and R.~N\"urnberg}, {\em A novel finite element
  approximation of anisotropic curve shortening flow}, Interfaces Free Bound.,
  25 (2023), pp.~671--708.

\bibitem{eqdproc}
\leavevmode\vrule height 2pt depth -1.6pt width 23pt, {\em An unconditionally
  stable finite element scheme for anisotropic curve shortening flow}, Arch.
  Math. (Brno), 59 (2023), pp.~263--274.

\bibitem{Dziuk94}
{\sc G.~Dziuk}, {\em Convergence of a semi-discrete scheme for the curve
  shortening flow}, Math. Models Methods Appl. Sci., 4 (1994), pp.~589--606.

\bibitem{Dziuk99}
\leavevmode\vrule height 2pt depth -1.6pt width 23pt, {\em Discrete anisotropic
  curve shortening flow}, SIAM J. Numer. Anal., 36 (1999), pp.~1808--1830.

\bibitem{EidelmanIK04}
{\sc S.~D. Eidelman, S.~D. Ivasyshen, and A.~N. Kochubei}, {\em Analytic
  methods in the theory of differential and pseudo-differential equations of
  parabolic type}, vol.~152 of Operator Theory: Advances and Applications,
  Birkh\"{a}user Verlag, Basel, 2004.

\bibitem{ElliottF17}
{\sc C.~M. Elliott and H.~Fritz}, {\em On approximations of the curve
  shortening flow and of the mean curvature flow based on the {D}e{T}urck
  trick}, IMA J. Numer. Anal., 37 (2017), pp.~543--603.

\bibitem{FonsecaM91}
{\sc I.~Fonseca and S.~M{\"u}ller}, {\em A uniqueness proof for the {W}ulff
  theorem}, Proc. Roy. Soc. Edinburgh Sect. A, 119 (1991), pp.~125--136.

\bibitem{Giga06}
{\sc Y.~Giga}, {\em Surface evolution equations}, vol.~99 of Monographs in
  Mathematics, Birkh\"{a}user, Basel, 2006.

\bibitem{Gurtin93}
{\sc M.~E. Gurtin}, {\em Thermomechanics of Evolving Phase Boundaries in the
  Plane}, Oxford Mathematical Monographs, The Clarendon Press Oxford University
  Press, New York, 1993.

\bibitem{Li20}
{\sc B.~Li}, {\em Convergence of {D}ziuk's linearly implicit parametric finite
  element method for curve shortening flow}, SIAM J. Numer. Anal., 58 (2020),
  pp.~2315--2333.

\bibitem{MaC07}
{\sc L.~Ma and D.~Chen}, {\em Curve shortening in a {R}iemannian manifold},
  Ann. Mat. Pura Appl. (4), 186 (2007), pp.~663--684.

\bibitem{Pozzi07}
{\sc P.~Pozzi}, {\em Anisotropic curve shortening flow in higher codimension},
  Math. Methods Appl. Sci., 30 (2007), pp.~1243--1281.

\bibitem{Pozzi12}
\leavevmode\vrule height 2pt depth -1.6pt width 23pt, {\em On the gradient flow
  for the anisotropic area functional}, Math. Nachr., 285 (2012), pp.~707--726.

\bibitem{Alberta}
{\sc A.~Schmidt and K.~G. Siebert}, {\em Design of Adaptive Finite Element
  Software: The Finite Element Toolbox {ALBERTA}}, vol.~42 of Lecture Notes in
  Computational Science and Engineering, Springer-Verlag, Berlin, 2005.

\bibitem{TaylorCH92}
{\sc J.~E. Taylor, J.~W. Cahn, and C.~A. Handwerker}, {\em Geometric models of
  crystal growth}, Acta Metall. Mater., 40 (1992), pp.~1443--1474.

\bibitem{Wheeler73}
{\sc M.~F. Wheeler}, {\em A priori {$L\sb{2}$} error estimates for {G}alerkin
  approximations to parabolic partial differential equations}, SIAM J. Numer.
  Anal., 10 (1973), pp.~723--759.

\bibitem{YeC21}
{\sc C.~Ye and J.~Cui}, {\em Convergence of {D}ziuk's fully discrete linearly
  implicit scheme for curve shortening flow}, SIAM J. Numer. Anal., 59 (2021),
  pp.~2823--2842.

\end{thebibliography}
\end{document}